\documentclass[secthm,seceqn,amsthm,ussrhead,10pt]{amsart}
\usepackage[utf8]{inputenc}
\usepackage[english]{babel}
\usepackage{amssymb,amsmath,amsthm,amsfonts,xcolor,enumerate,hyperref,comment,longtable,cleveref}

\usepackage{times}
\usepackage{cite}
\usepackage{pdflscape}
\usepackage{ulem}
\usepackage[mathcal]{euscript}
\usepackage{tikz}
\usepackage{hyperref}
\usepackage{cancel}
\usepackage{stmaryrd}
\usepackage[all,poly,graph]{xy}

\newcommand\Tstrut{\rule{0pt}{4.0ex}}       
\newcommand\Bstrut{\rule[-3.1ex]{0pt}{0pt}} 

\newcommand\bstrut{\rule[-1.5ex]{0pt}{0pt}} 

\usetikzlibrary{arrows}

\newtheorem{theorem}{Theorem}
\newtheorem{lemma}[theorem]{Lemma}

\newtheorem{remark}[theorem]{Remark}

\newtheorem{definition}[theorem]{Definition}
\newtheorem{corollary}[theorem]{Corollary}
\newtheorem{notation}[theorem]{Notation}

\newenvironment{Proof}[1][Proof.]{\begin{trivlist}
\item[\hskip \labelsep {\bfseries #1}]}{\flushright
$\Box$\end{trivlist}}

\usepackage{stmaryrd}
\usepackage{xcolor}

\begin{document}

{\LARGE  Classification of bilinear maps with   radical of codimension 2}\footnote{The work is  supported by the PCI of the UCA `Teor\'\i a de Lie y Teor\'\i a de Espacios de Banach', by the PAI with project numbers FQM298, FQM7156 and by the project of the Spanish Ministerio de Educaci\'on y Ciencia  MTM2016-76327C31P, RFBR 17-51-04004 and FAPESP 17/15437-6. }

   \

   {\bf Antonio Jesús Calder\'on$^{a}$, Amir Fern\'andez Ouaridi$^{a}$, Ivan Kaygorodov$^{b}$ \\

    \medskip
}

{\tiny

$^{a}$ Universidad de C\'adiz. Puerto Real, C\'adiz, Espa\~na.

$^{b}$ CMCC, Universidade Federal do ABC. Santo Andr\'e, Brasil.

\

\smallskip

   E-mail addresses:

\smallskip
    Antonio Jesús Calder\'on (ajesus.calderon@uca.es),

\smallskip
    Amir Fern\'andez Ouaridi (amir.fernandezouaridi@alum.uca.es),

\smallskip

    Ivan Kaygorodov (kaygorodov.ivan@gmail.com).

}

\

\

\noindent {\bf Abstract.}
Let ${\mathbb V}$ be an $n$-dimensional linear space over an algebraically closed base field. We provide a classification, up to equivalence, of all of the bilinear maps $f:{\mathbb V} \times {\mathbb V} \to {\mathbb V}$ such  that
${\rm dim}({\rm rad}(f)) =n-2$. This is equivalent to give  a complete classification (up to isomorphism) of all $n$-dimensional
algebras with annihilator of dimension  $n-2$ or, in other words,
a classification of  the  annihilator  extensions of all $2$-dimensional  algebras.

\section{Introduction}

 \begin{definition}\rm
 Let  $f:{\mathbb V} \times {\mathbb V} \to {\mathbb V}$ and $g:{\mathbb W} \times {\mathbb W}  \to {\mathbb W} $ be two   bilinear maps on the linear spaces  ${\mathbb V} $ and ${\mathbb W} $ (over the same base field) respectively.   We say that $f$ and $g$ are {\it equivalent} if there exists a linear isomorphism $\phi:{\mathbb V}  \to {\mathbb W} $ such that $\phi f(v_1,v_2)=g(\phi(v_1), \phi(v_2))$ for any $v_1, v_2 \in {\mathbb V} $. That is, if the following
 diagram is commutative

 $$\hbox{  \xymatrix{
            {\mathbb V} \times {\mathbb V} \ar[d]^{\phi \times \phi} \ar[r]^{\,\,\,\,\,\,f} &{\mathbb V} \ar[d]^{\phi}\\
   {\mathbb W} \times {\mathbb W} \ar[r]^{\,\,\,\,\,\,\,\,g}          &{\mathbb W}   .}}$$


 \end{definition}

 \begin{definition}\rm
Let ${\mathbb V}$ and ${\mathbb W} $ be  linear spaces (over a same base field) and   $f:{\mathbb V} \times {\mathbb V} \to {\mathbb W}$ a bilinear map. The {\it radical}  of $f$ is the set
$$ {\rm rad}(f)=\{ v\in {\mathbb V}: f(v,{\mathbb V})+ f({\mathbb V},v)=0  \}.$$
\end{definition}

\medskip

In the second  section of this work, and by inspiring in cohomology techniques in the study of Lie, Jordan and Malcev algebras (see \cite{ha17,hac16,ss78}), we develop a procedure to classify, up to equivalence, all of the bilinear maps
$f:{\mathbb V} \times {\mathbb V} \to {\mathbb V}$ with ${\rm dim}({\mathbb V}) =n$ and ${\rm dim}({\rm rad}(f))=p \neq 0$ in case the classification of all of the bilinear maps $f:{\mathbb W} \times {\mathbb W} \to {\mathbb W}$ when ${\rm dim}({\mathbb W}) =n-p$ was previously known.

\medskip

In  Section 3,  taking into account that the classification of the bilinear maps $f:{\mathbb W} \times {\mathbb W} \to {\mathbb W}$ when ${\rm dim}({\mathbb W}) =2$ is well-known, we apply the procedure  obtained in Section 2 to get a classification of all of the bilinear maps $f:{\mathbb V} \times {\mathbb V} \to {\mathbb V}$ when ${\rm dim}({\mathbb V}) =n$ and   ${\rm dim}({\rm rad}(f))=n-2$.

\medskip

Now we recall that an {\it algebra} $A$ is just a linear space, over a base field ${\bf k}$, endowed with a bilinear map, denoted by
$$\hbox{$\cdot:A \times A \to A$, $(x,y) \mapsto x \cdot y,$}$$ called the {\it product} of the algebra. By depending of the identities satisfied by the product we can speak about different categories of algebras (commutative, associative, Lie, Jordan, Malcev, alternative, and so on). The {\it dimension} of an algebra is its dimension as a linear space and the {\it annihilator} of $A$ is defined as the set $${\rm ann}(A)= \{x \in A: x \cdot A + A \cdot x=0\}.$$ Taking into account this observation we have that  the present paper gives us  a complete classification of all of the $n$-dimensional
algebras, over an algebraically closed field, with annihilator of dimension  $n-2$. This is  the classification of all of the annihilator extensions of an arbitrary 2-dimensional algebra. Here we recall that annihilator extensions of algebras, called central extension in the (anti)commutative case, play an important role in quantum mechanics: one of the earlier
encounters is by means of Wigner’s theorem which states that a symmetry of a quantum
mechanical system determines an (anti-)unitary transformation of a Hilbert space.
Central extensions are needed in physics, because the symmetry group of a quantized
system usually is a central extension of the classical symmetry group, and in the same way
the corresponding symmetry Lie algebra of the quantum system is, in general, a central
extension of the classical symmetry algebra.
The Virasoro
algebra is the universal central extension of the Witt algebra, the Heisenberg algebra is
the central extension of a commutative Lie algebra  \cite[Chapter 18]{bkk}.
The algebraic study of central extensions of Lie and non-Lie algebras has a very big story \cite{omirov,zusmanovich,ha17,hac16,is11,ss78}.
So, Skjelbred and Sund used central extensions of Lie algebras for a classification of nilpotent Lie algebras  \cite{ss78}.
After that, using the method described by Skjelbred and Sund were classified
different kind of algebras \cite{hegazi1,hegazi2,ha17,hac16}.

\medskip

Throughout all of the paper we will consider bilinear maps  $f:{\mathbb V} \times {\mathbb V} \to {\mathbb V}$  where ${\mathbb V}$ is an $n$-dimensional linear space over an algebraically closed base field ${\bf k}$. For a more comfortable development of our study we introduce the next concept:

\begin{definition}\rm
A {\it bilinear pair} of dimension $n$   is a pair $({\mathbb V},f)$ in which ${\mathbb V}$ is an $n$-dimensional linear space and $f:{\mathbb V} \times {\mathbb V} \to {\mathbb V}$ a bilinear map. The {\it radical}  of $({\mathbb V},f)$ is the radical of $f$.
\end{definition}

\medskip

We recall that in the theory of Lie algebras, an {\it annihilator component} $U$ of a Lie algebra $L$ is an one-dimensional ideal of $L$ such that there exists another ideal
  $W$ of $L$ satisfying $L=W \oplus U$. Observe that the product in $L$ is then $[(w_1,u_1),(w_2,u_2)]=([w_1,w_2],0)$ for any $w_1,w_2 \in W$ and $u_1,u_2 \in U$. This concept can be extended in our framework as follows:

    \begin{definition}\rm
A {\it radical component} of a bilinear pair $({\mathbb V},f)$ is an one-dimensional linear subspace ${\mathbb U}$ of ${\mathbb V}$ such that $f({\mathbb U},{\mathbb V})+f({\mathbb V},{\mathbb U})=0$  and  satisfying that there exists  a linear subspace ${\mathbb W}$ of ${\mathbb V}$ with
${\mathbb V}={\mathbb W} \oplus {\mathbb U}$   and   $f({\mathbb W},{\mathbb W}) \subset {\mathbb W}$. We will write $({\mathbb V},f)=({\mathbb W} \oplus {\bf k}, f|_{{\mathbb W} \times {\mathbb W}})$.
\end{definition}


\begin{remark}\rm
Observe that if ${\mathbb U}$ is a radical  component of a bilinear pair $({\mathbb V},f)$, then by choosing some   $0 \neq u \in {\mathbb U}$, by identifying  ${\mathbb U}={\bf k}u:= {\bf k}$ and by writing ${\mathbb V}:={\mathbb W} \oplus {\bf k}$ we get $f(w_1+\lambda, w_2+ \beta)=f(w_1,w_2)\in {\mathbb W}$.  This reduces the study of $({\mathbb V},f)$ to those of $({\mathbb W},f|_{{\mathbb W} \times {\mathbb W}})$ and allows us to write $({\mathbb V},f)=({\mathbb W} \oplus {\bf k}, f|_{{\mathbb W} \times {\mathbb W}}).$
\end{remark}



The main aim of the present paper is to prove the following theorem:

\begin{theorem}[Main Theorem]\label{maintheo}
Let $({\mathbb V},f)$ be an $n$-dimensional bilineal pair  with an $(n-2)$-dimensional radical.

\begin{itemize}
\item \quad If $n=3$, then $({\mathbb V},f)$ is equivalent to one of the following non-equivalent bilinear pairs. Either $({\mathbb V},f)$ is equivalent to ${({\bf A}' \oplus {\bf k}, f)}$, where $({\bf A}', f|_{{\bf A}' \times {\bf A}'})$ is a $2$-dimensional bilinear pair (see \hyperref[tab4]{Table 4} and Remark \ref{remrad}) with zero radical; or to $({\bf A}_i=({\mathbb V},f_i))$ for some $i=1, \ldots, 68,$ where any $f_i$ is defined  (for a fixed   basis ${\mathcal B}=\{e_1,e_2,e_3\}$ of  ${\mathbb V}$) as    $f_i(e_j,e_k):=e_je_k$ in \hyperref[tab1]{Table 1}.
\medskip
\item \quad If $n=4,$ then $({\mathbb V},f)$ is equivalent to  one of the following non-equivalent bilinear pairs.  Either $({\mathbb V},f)$ is equivalent to ${({\bf A}' \oplus {\bf k}\oplus {\bf k}, f)}$, where $({\bf A}', f|_{{\bf A}' \times {\bf A}'})$ is a $2$-dimensional bilinear pair (see \hyperref[tab4]{Table 4}  and Remark \ref{remrad}) with zero radical; or to $({\bf A}_i \oplus {\bf k}, f_i)$ for some $i=1, \ldots, 68;$
or to $({\bf A}_i=({\mathbb V},f_i))$ for some $i=69, \ldots, 121,$  where any $f_i$ is defined  (for a fixed   basis ${\mathcal B}=\{e_1,e_2,e_3,e_4\}$ of ${\mathbb V}$) as    $f_i(e_j,e_k):=e_je_k$ in \hyperref[tab2]{Table 2}.
\medskip
\item \quad If $n=5,$  then $({\mathbb V},f)$ is equivalent to  one of the following non-equivalent bilinear pairs.  Either $({\mathbb V},f)$ is equivalent to ${({\bf A}' \oplus {\bf k}\oplus {\bf k}\oplus {\bf k}, f)}$, where $({\bf A}', f|_{{\bf A}' \times {\bf A}'})$ is a $2$-dimensional bilinear pair (see \hyperref[tab4]{Table 4} and Remark \ref{remrad}) with zero radical; or to $({\bf A}_i \oplus {\bf k}\oplus {\bf k}, f_i)$ for some $i=1, \ldots, 68$;
or to $({\bf A}_i\oplus {\bf k},f_i)$ for some $i=69, \ldots, 121;$ or to  $({\bf A}_i=({\mathbb V},f_i))$ for some $i=122, \ldots, 133$ where any $f_i$ is defined  (for a fixed   basis ${\mathcal B}=\{e_1,e_2,e_3,e_4,e_5\}$ of ${\mathbb V}$) as    $f_i(e_j,e_k):=e_je_k$  in \hyperref[tab3]{Table 3}.
\medskip
\item \quad If $n=6,$   then $({\mathbb V},f)$ is equivalent to  one of the following non-equivalent bilinear pairs.  Either $({\mathbb V},f)$ is equivalent to ${({\bf A}' \oplus {\bf k}\oplus {\bf k}\oplus {\bf k}\oplus {\bf k}, f)}$, where $({\bf A}', f|_{{\bf A}' \times {\bf A}'})$ is a $2$-dimensional bilinear pair (see \hyperref[tab4]{Table 4} and Remark \ref{remrad}) with zero radical; or to $({\bf A}_i \oplus {\bf k}\oplus {\bf k}\oplus {\bf k}, f_i)$ for some $i=1, \ldots, 68$;
or to $({\bf A}_i\oplus {\bf k}\oplus {\bf k},f_i)$ for some $i=69, \ldots, 121;$ or to  $({\bf A}_i\oplus {\bf k},f_i)$ for some $i=122, \ldots, 133;$ or to
$({\bf A}_{134}=({\mathbb V},f_{134}))$ where  $f_{134}$ is defined  (for a fixed   basis ${\mathcal B}=\{e_1,e_2,e_3,e_4,e_5, e_6\}$ of ${\mathbb V}$) as:
$$\begin{array}{lllll lll}
{\bf A}_{134}&:&(\mathfrak{N}_2)_{6,14} &:&    e_1 e_1 = e_3 & e_1 e_2= e_4 & e_2 e_1=e_5  & e_2 e_2=e_6
\end{array}$$
\item \quad If $n \geq 7$ then $({\mathbb V},f)$ has at least $(n-6)$ radical components.
\end{itemize}

\end{theorem}

\section{Development of the techniques}


Throughout the paper, given two linear spaces $\mathbb U$ and $\mathbb W$ over a same base field, we will denote by $Hom(\mathbb U, \mathbb W)$ the set of all of the linear maps from $\mathbb U$ to $ \mathbb W$ and by  $Bil(\mathbb U\times \mathbb U, \mathbb W)$ the set of all of the bilinear maps from $\mathbb U\times \mathbb U$ onto $\mathbb W$.

\medskip

\begin{notation}\label{pig1}
Let $(\mathbb U,g)$ be
 a bilinear pair,  $\mathbb W$ a linear space over the same base field and $\theta \in Bil(\mathbb U\times \mathbb U, \mathbb W)$.
We denote  $(\mathbb U\oplus \mathbb W, g_{\theta})$ the bilinear pair  given by
 \begin{align*}
g_{\theta}: (\mathbb U\oplus \mathbb W) \times (\mathbb U\oplus \mathbb W) &\rightarrow \mathbb U \oplus \mathbb W \\
((u_1,w_1),(u_2,w_2)) &\mapsto (g(u_1,u_2), \theta(u_1, u_2)).
\end{align*}

\medskip

Observe that ${\rm rad}(g_{\theta})=({\rm rad}(g)\cap {\rm rad}(\theta))\oplus \mathbb W$.
\end{notation}

\begin{definition}
\label{2.4}Given a bilinear pair $(\mathbb U,g)$ and a linear space $\mathbb W$. For any $h\in Hom(\mathbb U, \mathbb W)$ we define $\delta  h \in Bil(\mathbb U\times \mathbb U, \mathbb W)$ as
\begin{align*}
\delta  h :\mathbb U\times \mathbb U & \rightarrow \mathbb W  \\
(u_1,u_2)&\mapsto h(g(u_1,u_2)).
\end{align*}

We  also introduce  the linear
subspace $$B_g(\mathbb U\times \mathbb U, \mathbb W)=
\{ \delta  h \in Bil(\mathbb U\times \mathbb U, \mathbb W): h \in Hom(\mathbb U, \mathbb W)\}.$$

\end{definition}

\begin{notation}
Given a   bilinear map $f:\mathbb V \times \mathbb V \to \mathbb V$ we will denote by $\tilde{f}:({\mathbb V} / {{\rm rad}(f)}) \times ({\mathbb V} / {{\rm rad}(f)})\to {\mathbb V} / {{\rm rad}(f)}$ the (well-defined) bilinear map given by $\tilde{f}([x],[y]):=[f(x,y)]$ for any $[x],[y] \in {\mathbb V} / {{\rm rad}(f)}$.
\end{notation}

\begin{lemma}\label{l22}
Let $(\mathbb V,f)$ be a bilinear pair of dimension $n$ such that
${\rm dim}({\rm rad}(f))=m\not= 0$. Then there exist
 a unique (up to equivalence) bilinear
pair $(\mathbb U,g)$ of dimension $n-m$, and a bilinear map
$\theta:\mathbb U\times \mathbb U \rightarrow {\rm rad}(f)$ with   ${\rm rad}(g)\cap {\rm rad}(\theta)=0$,
in such a way that
$(\mathbb V,f)$ is equivalent to $(\mathbb U\oplus {\rm rad}(f), g_{\theta})$ and that $({\mathbb V} / {{\rm rad}(f)}, \tilde{f})$ is equivalent  to $(\mathbb U,g)$.
\end{lemma}

\begin{Proof}
Consider a linear complement $\mathbb U$ of ${\rm rad}(f)$ in $\mathbb V$, being so
 $\mathbb V=\mathbb U\oplus {\rm rad}(f)$, and define
  $$g=\pi_{\mathbb U}\circ f: \mathbb U \times \mathbb U \rightarrow \mathbb U$$ where
$\pi_{\mathbb U}$ is the projection onto $\mathbb U$, and
$$\theta=\pi_{{\rm rad}(f)}\circ f: \mathbb U \times \mathbb U \rightarrow {\rm rad}(f)$$ where $\pi_{{\rm rad}(f)}$ is the
projection onto ${\rm rad}(f)$. Observe that
 ${\rm rad}(g)\cap {\rm rad}(\theta)=0$.

 \smallskip

Since for any   $x,y \in \mathbb V$ we have   $$(\pi_{\mathbb U} f(x,y),\pi_{{\rm rad}(f)} f(x,y))=$$
\begin{equation}\label{pre1}
(\pi_{\mathbb U} f(\pi_{\mathbb U}(x)+\pi_{{\rm rad}(f)}(x),\pi_{\mathbb U}(y)+\pi_{{\rm rad}(f)}(y)),
\pi_{{\rm rad}(f)} f(\pi_{\mathbb U}(x)+\pi_{{\rm rad}(f)}(x),\pi_{\mathbb U}(y)+\pi_{{\rm rad}(f)}(y))=
\end{equation}
 $$(\pi_{\mathbb U} f(\pi_{\mathbb U}(x),\pi_{\mathbb U}(y)),\pi_{{\rm rad}(f)} f(\pi_{\mathbb U}(x),\pi_{\mathbb U}(y))),$$

then the next
 diagram is commutative
 $$\hbox{  \xymatrix{
            V \times V \ar[d]^{(\pi_{\mathbb U},\pi_{{\rm rad}(f)}) \times (\pi_{\mathbb U},\pi_{{\rm rad}(f)})} \ar[r]^{f} &V \ar[d]^{(\pi_{\mathbb U},\pi_{{\rm rad}(f)})}\\
   (\mathbb U\oplus {\rm rad}(f)) \times (\mathbb U\oplus {\rm rad}(f))    \ar[r]^{ \,\,\,\,\,\,\,\,\,\,\,\,\,\,\,\,\,\,\,\,\,\,\,\,\,\,\,\,  g_{\theta}}          &\mathbb U\oplus {\rm rad}(f)   }}$$

   \medskip

   and so the bilinear pairs$(\mathbb V,f)$ and $(\mathbb U\oplus {\rm rad}(f), g_{\theta})$ are  equivalent.

   \medskip

   Denote by $\tilde{\pi_{\mathbb U}}:{\mathbb V}/{{\rm rad}(f)}  \to {\mathbb U}$ the linear map given by $\tilde{\pi_{\mathbb U}}([x]):=\pi_{\mathbb U}(x)$ for any $[x] \in {\mathbb V}/{{\rm rad}(f)}$. Then Equation (\ref{pre1}) gives us that the next diagram is commutative:
   $$\hbox{  \xymatrix{
            ({\mathbb V}/{{\rm rad}(f)}) \times ({\mathbb V}/{{\rm rad}(f)}) \ar[d]^{\tilde{\pi_{\mathbb U} }\times \tilde{ \pi_{\mathbb U}}} \ar[r]^{\,\,\,\,\,\,\,\,\,\,\,\,\,\,\,\,\,\,\,\,\,\,\,\,\,\,\tilde{f}} & {\mathbb V}/{{\rm rad}(f)} \ar[d]^{\tilde{ \pi_{\mathbb U}}}\\
  { \mathbb U} \times {\mathbb U}    \ar[r]^{g}          &{\mathbb U}   }}$$

   \medskip

  From here,  the bilinear pairs $({\mathbb V} / {{\rm rad}(f)}, \tilde{f})$ and $(\mathbb U,g)$ are equivalent and we also have the uniqueness of $(\mathbb U,g)$.

\end{Proof}

\begin{remark}\label{pig2}
At this point, observe from Notation \ref{pig1} and Lemma \ref{l22}, that
  $(\mathbb V,f)$ is a bilinear pair of dimension $n$ with
${\rm dim}({\rm rad}(f))=m\not= 0$ if and only if   there exist
 a  bilinear
pair $(\mathbb U,g)$ of dimension $n-m$ and a bilinear map
$\theta:\mathbb U\times \mathbb U \rightarrow \mathbb W$, where $\mathbb W$ is a linear space of dimension $m$, with   ${\rm rad}(g)\cap {\rm rad}(\theta)=0$,
in such a way that
$(\mathbb V,f)$ is equivalent to $(\mathbb U\oplus \mathbb W, g_{\theta})$.

\smallskip

From here, in case we knew a classification  of the  bilinear
pairs $(\mathbb U,g)$ of dimension $n-m$. The problem of classifying the bilinear pairs of dimension $n$ with
${\rm dim}({\rm rad}(f))=m$ would be  reduced to study the set $\{\theta \in Bil(\mathbb U\times \mathbb U,\mathbb W) :{\rm rad}(g)\cap {\rm rad}(\theta)=0\}$ where  $\mathbb W$ is a linear space of ${\rm dim}(\mathbb W)=m$. However,  Lemma \ref{l23} will allow us to simplify our study.
\end{remark}

\begin{lemma}
\label{l23}Let $(\mathbb U,g)$ be a bilinear pair,  $\mathbb W$  a vector space and
$\theta, \mu : \mathbb U\times\mathbb  U \rightarrow \mathbb W$ two bilinear maps such that
$\theta-\mu \in B_g(\mathbb U \times {\mathbb U },\mathbb W)$. Then the bilinear pairs
$(\mathbb U\oplus \mathbb W, g_{\theta})$ and $(\mathbb U\oplus \mathbb W, g_{\mu})$ are equivalent and
${\rm rad}(g_{\theta})={\rm rad}(g_{\mu})$.
\end{lemma}

\begin{Proof}
We know there exists a linear map $h:\mathbb U\rightarrow \mathbb W$ such that
$\theta- \mu=\delta h$. Let us define   the linear isomorphism
\begin{align*}
\sigma:\mathbb U\oplus \mathbb W &\rightarrow \mathbb U\oplus\mathbb  W  \\
(u,w)&\mapsto  (u,h(u)+w).
\end{align*}

Then we have
$\sigma(g_{\mu}(u_1+w_1,u_2+w_2))=\sigma(g(u_1,u_2),\mu(u_1,u_2))=(g(u_1,u_2),h(g(u_1,u_2))+\mu(u_1,u_2))$ and also
$g_{\theta}((u_1,h(u_1)+w_1),(u_2,h(u_2)+w_2))=(g(u_1,u_2),\theta(u_1,u_2))$ for any $u_1,u_2\in \mathbb U$ and $w_1,w_2\in \mathbb W$.
From here, $\sigma g_{\mu}= g_{\theta} (\sigma \times \sigma)$ and so the bilinear pairs $(\mathbb U\oplus \mathbb W, g_{\theta})$ and $(\mathbb U\oplus \mathbb W, g_{\mu})$ are equivalent.

\smallskip

Finally, observe that $\theta(u_1,u_2)=\mu(u_1,u_2)+h(g(u_1,u_2))$ for any $u_1,u_2\in \mathbb U$.
 Thus, $\mu(u_1,u_2)=g(u_1,u_2)=0$ if and only if
$\theta(u_1,u_2)=g(u_1,u_2)=0$  and so  ${\rm rad}(g_{\mu})={\rm rad}(g_{\theta})$.
\end{Proof}

This result leads us to introduce the next concept:

\begin{definition}
For any bilinear pair $(\mathbb U,g)$ and any  linear space $\mathbb W$, we call {\it second cohomology space} of $(\mathbb U,g)$ respect to   $\mathbb W$,  to the linear space $$H^2_g(\mathbb U,\mathbb W):= \cfrac{Bil(\mathbb U\times \mathbb U,\mathbb W)}{B_g(\mathbb U \times {\mathbb U },\mathbb W)}.$$
\end{definition}

\begin{remark}\label{pig2}
Now we observe that Remark \ref{pig2} and Lemma \ref{l23} give us that   in case we knew a classification  of the  bilinear
pairs $(\mathbb U,g)$ of dimension $n-m$. The problem of classifying the bilinear pairs of dimension $n$ with
${\rm dim}({\rm rad}(f))=m$ would be  reduced to study the sets $\{[\theta] \in H^2_g(\mathbb U,\mathbb W) :{\rm rad}(g)\cap {\rm rad}(\theta)=0\}$ where  $\mathbb W$ is a linear space of ${\rm dim}(\mathbb W)=m$. However,  we have not  enough information to distinguish when two different elements $[\theta]$ and $[\mu]$ in the above set give rise to non-equivalent bilinear pairs. Lemma \ref{l24} will clarify this point.
\end{remark}

\begin{notation}
Let $\mathbb U$ and  $\mathbb W$ be a couple of  linear spaces, $\mu \in Bil(\mathbb U \times \mathbb U, \mathbb W)$ and
$\phi \in Hom(\mathbb U , \mathbb U)$, then we will denote by $\phi\mu \in Bil(\mathbb U \times \mathbb U, \mathbb W)$ the bilinear map given by
$\phi\mu(u_1,u_2):=\mu(\phi(u_1), \phi(u_2))$ for any $u_1,u_2 \in \mathbb U$.
\end{notation}

\begin{lemma}
\label{l24}
Let $(\mathbb U,g)$ be a bilinear pair,  $\mathbb W$  a vector space and
$\theta, \mu : \mathbb U\times\mathbb  U \rightarrow \mathbb W$ two bilinear maps such that ${\rm rad}(g)\cap {\rm rad}(\theta)={\rm rad}(g)\cap {\rm rad}(\mu)=0$. Then we have
that the bilinear pairs   $(\mathbb U\oplus \mathbb W, g_{\theta})$ and $ (\mathbb U\oplus \mathbb W, g_{\mu})$ are equivalent  if and only if
there exist a linear isomorphism $\phi:\mathbb U\rightarrow \mathbb U$ satisfying $\phi g= g (\phi \times \phi)$, and a linear
isomorphism $\psi: \mathbb W\rightarrow \mathbb W$ such that
$\phi\mu - \psi \theta \in B_g(\mathbb U \times {\mathbb U }, \mathbb W)$.

\end{lemma}

\begin{Proof}
Suppose the bilinear pairs
$(\mathbb U\oplus \mathbb W, g_{\theta})$ and $(\mathbb U\oplus \mathbb W, g_{\mu})$ are equivalent. Then there exists a linear isomorphism  $\Phi$
such that $\Phi g_{\theta}= g_{\mu} (\Phi \times \Phi)$.
Since $\Phi(\mathbb W)=\Phi({\rm rad}(g_{\theta}))={\rm rad}(g_{\mu})=\mathbb W$ then
$\psi:={{  \left.\kern-\nulldelimiterspace   \phi     \right|_{\mathbb W}   }}$
is a linear isomorphism.
Let $\{e_i\}_{i=1, \ldots ,n}$ be a basis of $\mathbb U$ and
write  $\Phi(e_i)=e'_i+w_i$, where $e'_i\in\mathbb  U$ and $w_i\in \mathbb W$. Now
define the linear isomorphisms:
\begin{align*}
\phi :\mathbb  U &\rightarrow\mathbb  U &\varphi : \mathbb U &\rightarrow\mathbb  W \\
e_i &\mapsto e'_i. & e_i &\mapsto w_i
\end{align*}


For $u_1, u_2 \in \mathbb U$ we have
$$\Phi(g_{\theta}(u_1, u_2))=\Phi(g(u_1,u_2)+\theta(u_1,u_2))= \phi(g(u_1,u_2))+ \varphi(g(u_1,u_2))+\psi(\theta(u_1,u_2))$$
and
$$g_{\mu}(\Phi(u_1),\Phi(u_2))=g_{\mu}(\phi(u_1)+\varphi(u_1),\phi(u_2)+\varphi(u_2))=g(\phi(u_1),\phi(u_2))+\mu(\phi(u_1),\phi(u_2)).$$

Since $\Phi g_{\theta}= g_{\mu} (\Phi \times \Phi)$ we can match both equations. The part in
$\mathbb U$ shows that $\phi$ satisfies $\phi g= g (\phi \times \phi)$,  and the part in $\mathbb W$ shows that
$$\mu(\phi(u_1),\phi(u_2))=\varphi(g(u_1,u_2))+\psi(\theta(u_1,u_2)).$$

Thus, $\phi\mu - \psi \theta =\delta \varphi \in B_g(\mathbb U \times {\mathbb U },\mathbb W)$.

\medskip

Conversely,  in case there exist a linear isomorphism $\phi:\mathbb U\rightarrow \mathbb U$ satisfying $\phi g= g (\phi \times \phi)$, and a linear
isomorphism $\psi: \mathbb W\rightarrow \mathbb W$ such that
$\phi\mu - \psi \theta \in B_g(\mathbb U \times {\mathbb U }, \mathbb W)$. By writing  $\phi\mu - \psi \theta = \delta \varphi$ for $\varphi \in Hom({\mathbb U }, \mathbb W)$ we
 can define the linear isomorphism $\Phi: \mathbb U\oplus \mathbb W \to \mathbb U\oplus \mathbb W$ as
$$\Phi:= \left( {\begin{array}{cc}
   \phi & 0 \\
   \varphi & \psi
  \end{array} } \right)$$

which satisfies  $ \Phi g_{\theta}=g_{\mu} (\Phi \times \Phi)$ and so the bilinear pairs $(\mathbb U\oplus \mathbb W, g_{\theta})$ and $(\mathbb U\oplus \mathbb W, g_{\mu})$ are equivalent.
\end{Proof}

\begin{remark}\label{pig3}
Taking into account Remark \ref{pig2} and Lemma \ref{l24},  we have a first procedure to classify (up to equivalence)  the bilinear pairs of dimension $n$ with
${\rm dim}({\rm rad}(f))=m \neq 0$, in case we  knew the classification  of the  bilinear
pairs $(\mathbb U,g)$ of dimension $n-m$. This will follow three steps for any pair $(\mathbb U,g)$:

\begin{itemize}
\item[(i)] Compute the set \hbox{${\mathcal A}:=\{\phi \in Hom(\mathbb U,\mathbb U): \phi$ is bijective and $\phi g=g(\phi \times \phi)\}$}.
\item[(ii)] Compute $H^2_g(\mathbb U, \mathbb W)$ where ${\rm dim}(\mathbb W)=m$.
\item[(iii)] For any $\phi \in {\mathcal A}$ and any linear isomorphism $\psi$ of  $\mathbb W$, study when $[\psi \mu]=[\phi \theta]$ in $H^2_g(\mathbb U, \mathbb W)$ where $\theta, \mu \in Bil(\mathbb U \times \mathbb U, \mathbb W )$ are such that ${\rm rad}(g)\cap {\rm rad}(\theta)={\rm rad}(g)\cap {\rm rad}(\mu)=0.$
\end{itemize}

 \medskip

 However, we can improve so much this procedure by introducing ``coordinates'' in $H^2_g(\mathbb U, \mathbb W)$ as follows.
\end{remark}



\medskip

Consider a    bilinear
pair $(\mathbb U,g)$, a linear space $\mathbb W$ and fix   a basis $\{e_{1},\ldots ,e_{m}\}$ of ${\mathbb W}$. For any $f\in Bil(\mathbb U\times \mathbb U, {\mathbb W})$ and  $u_1, u_2 \in \mathbb U$ we can write  $$f\left(u_1,u_2\right) =\underset{i=1}{\overset{m}{\sum }} f_{i}\left( u_1,u_2\right) e_{i}$$  where $f_{i}\in Bil(\mathbb U\times \mathbb U, {\bf k}) $ for $i \in \{1,...,m\}$. Observe that $${\rm rad}(f)={\rm rad}(f_1)\cap\cdots\cap {\rm rad}(f_m)$$
and that
$f\in B_g(\mathbb U \times \mathbb U, \mathbb W)$ if and only if ${f}_i\in B_{g}(\mathbb U \times \mathbb U,{\bf k})$ for any  $i \in \{1,...,m\}$.

\begin{lemma}
\label{l28}Let $(\mathbb U,g)$ be a bilinear pair and $\{u_1,\ldots, u_m\}$ a basis of
$g(\mathbb U,\mathbb U)$ (linear subspace generated by $\{g(x,y): x,y \in \mathbb U\}$). Then $\{\delta u_1^*,\ldots, \delta u_m^*\}$ is a basis
of $B_g(\mathbb U \times \mathbb U, {\bf k})$.
\end{lemma}

\begin{Proof}
By extending  $\{u_1,\ldots, u_m\}$ to a basis $\{u_i\}_{i\leq n}$ of $\mathbb U$ we get
the  basis $\{u_i^*\}_{i\leq n}$ of
$Hom(\mathbb U,{\bf k})$. Then for any  $\delta h \in B_g(\mathbb U \times \mathbb U, {\bf k})$ we can write
$h=\sum_{i=1}^{n} \alpha_i u_i^*$ for some $\alpha_i \in {\bf k}$ and so
$$\delta h(u_l, u_k)= \sum_{i=1}^{n} \alpha_i u_i^* g(u_l, u_k)= \sum_{i=1}^{n}\alpha_i u_i^* (\sum_{j=1}^{m}\beta_j u_j)= \sum_{i=1}^{m}\alpha_i u_i^* g(u_l,u_k)
= \sum_{i=1}^{m}\alpha_i \delta u_i^*(u_l,u_k).$$

Thus, $span(\{\delta u_1^*,\ldots, \delta u_m^*\})=B_g(\mathbb U\times \mathbb U,{\bf k})$. Let
$\alpha_1, \ldots, \alpha_m\in {\bf k}$ be such that
$\sum_{i=1}^{m}\alpha_i \delta u_i^*=0$, then
$$\sum_{i=1}^{m} \alpha_i \delta u_i^*(\mathbb U,\mathbb U)=\sum_{i=1}^{m}\alpha_i u_i^*(g(\mathbb U,\mathbb U))=0.$$

Since $\{u_i^*\}_{i\leq m}$ is a basis of $Hom(f(\mathbb U,\mathbb U), {\bf k})$ then
$\alpha_1=\cdots=\alpha_m=0$. Thus, we conclude
$\{\delta u_1^*,\ldots, \delta u_m^*\}$ is a basis of $B_g(\mathbb U\times \mathbb U, {\bf k})$.
\end{Proof}

\begin{corollary}
Let $(\mathbb U,g)$ be a bilinear pair, then ${\rm dim}(H^2(\mathbb U, {\bf k}))={\rm dim}(\mathbb U)^2-{\rm dim}(f(\mathbb U, \mathbb U))$.
\end{corollary}

\begin{lemma}
\label{l210}
Let $(\mathbb U,g)$ be a bilinear pair,
 $\mathbb W$ a linear space with a fixed basis $\{e_1,\ldots, e_m \}$, and  $\theta \in Bil(\mathbb U\times \mathbb U, \mathbb W)$ such that
$\theta(u_1, u_2)=\sum_{i=1}^{m}\theta_i(u_1,u_2)e_i$. Suppose
${\rm rad}(\theta)\cap {\rm rad}(g)=0$, then the bilinear pair $(\mathbb U\oplus \mathbb W, g_{\theta})$ has a radical
component if and only if the set $\{ \left[ \theta_i \right] \}_{i=1}^m$ is
linearly dependent in $H^2(\mathbb U, {\bf k})$.
\end{lemma}

\begin{Proof}
Suppose $(\mathbb U\oplus \mathbb W, g_{\theta})$ has a radical component. Then there exists a nonzero element $v_1 \in \mathbb W$ such that $g_{\theta}(v_1,\mathbb U\oplus \mathbb W)+g_{\theta}(\mathbb U\oplus \mathbb W,v_1)=0$ and
$\mathbb U\oplus \mathbb W= I\oplus \langle v_1 \rangle$ with    $g_{\theta}(I,I) \subset I$.

Enlarge the set $\{v_1\}$ to the set   $\{v_1,\ldots, v_m\}$ to form a basis of $\mathbb W$. Then  there exists an
invertible matrix $(a_{ij})_{ij}$ such that
$e_i=\sum_{j=1}^{m}a_{ij}v_j$ for any $i=1,...,m$.

Thus,
$\theta(u_1,u_2)=\sum_{i=1}^{m}\theta_i(u_1,u_2)(\sum_{j=1}^{m}a_{ij}v_j)=\sum_{j=1}^{m}(\sum_{i=1}^{m}\theta_i(u_1,u_2)a_{ij})v_j$.
Since $g_{\theta}(\mathbb U\oplus \mathbb W, \mathbb U\oplus \mathbb W)\subset I$ we have
$\sum_{i=1}^{m}a_{i,1}\theta_i(u_1,u_2)=0$ for all $u_1, u_2\in \mathbb U$.
Then $\sum_{i=1}^{m}a_{i1}\theta_i=0$ and therefore
$\sum_{i=1}^{m}a_{i1}\left[\theta_i\right]=0$. Now, if the set
 $\{\left[\theta_i\right]\}_{i\leq m}$ was  linearly independent,  then
$a_{11}=\cdots=a_{n1}=0$ and  $(a_{ij})_{ij}$ would not be invertible. So the set $\{\left[\theta_i\right]\}_{i\leq m}$ is linearly  dependent.

Conversely, suppose that the set  $\{\left[\theta_i\right]\}_{i\leq m}$ is  linearly
dependent in $H^2(\mathbb U, {\bf k})$. Without loss of generality we may assume that
$\left[\phi_m \right]=\sum_{i=1}^{m}\alpha_i\left[\theta_i \right]$ for
$\alpha_i \in {\bf k}$. Define
 $v(u_1,u_2)=\sum_{i=1}^{m} v_i(u_1,u_2)e_i$ by setting  $v_i=\theta_i$
for  $i\leq m-1$ and $v_m=\sum_{i=1}^{m-1}\alpha_i \theta_i$. Then
$\left[v \right]=\left[\theta \right]$ and so the bilinear pairs $(\mathbb U\oplus \mathbb W, g_{\theta})$ and $(\mathbb U\oplus \mathbb W, g_{v})$ are equivalent.
Also:
$$v(u_1,u_2)=\sum_{i=1}^{m-1} \theta_i(u_1,u_2)(e_i+\alpha_i e_m).$$

For $i=1,\ldots,m-1$ choose $w_i=e_i+\alpha_i e_m$, then
$v(u_1,u_2)=\sum_{i=1}^{m-1}\theta_i (u_1,u_2)w_i$.

Thus,
$g_v(\mathbb U\oplus \mathbb W, \mathbb U\oplus \mathbb W)\subset \mathbb U\oplus \langle w_1,\ldots,w_{m-1} \rangle$.
Then $(\mathbb U\oplus \mathbb W,g_{\theta})$ has a radical component.
\end{Proof}

\begin{lemma}
\label{l211}
Let $(\mathbb U,g)$ be a bilinear pair,
 $\mathbb W$ a linear space with a fixed basis $\{e_1,\ldots, e_m \}$, and  $\theta, \mu \in Bil(\mathbb U\times \mathbb U, \mathbb W)$ such that
$\theta(u_1, u_2)=\sum_{i=1}^{m}\theta_i(u_1,u_2)e_i$ and $\mu(u_1,u_2)= \sum_{i=1}^{m} \mu_i (u_1,u_2)e_i$. Suppose that
${\rm rad}(\theta)\cap {\rm rad}(g)={\rm rad}(\mu)\cap {\rm rad}(g)=0$ and that $(\mathbb U\oplus \mathbb W, g_{\theta})$ has not radical components. Then the bilinear pairs $(\mathbb U\oplus \mathbb W, g_{\theta})$
and $(\mathbb U\oplus \mathbb W, g_{\mu})$ are equivalent if and only if there
exists an isomorphism $\phi: \mathbb U\rightarrow \mathbb U$ such that the set
$\{\left[ \phi\theta_i \right]\}_{i\leq m}$ spans the same subspace of
$H^2(\mathbb U, {\bf k})$ than the set $\{\left[ \mu_i \right]\}_{i\leq m}$.
\end{lemma}

\begin{Proof}
Suppose  that the bilinear pairs
$(\mathbb U\oplus \mathbb W, g_{\theta})$ and $(\mathbb U\oplus \mathbb W, g_{\mu})$ are equivalent. By Lemma \ref{l24}, there exist  a linear isomorphism $\phi:\mathbb U\rightarrow \mathbb U$ satisfying $\phi g= g (\phi \times \phi)$, and a linear
isomorphism $\psi: \mathbb W\rightarrow \mathbb W$ such that
$\phi\mu - \psi \theta \in B_g(\mathbb U \times {\mathbb U }, \mathbb W)$.

Let
$\psi(e_i)=\sum_{j=1}^{m}a_{ij}e_j$. Then
$$(\phi \mu - \psi \theta)(u_1,u_2)=\sum_{j=1}^{m}(\phi\mu_j-\sum_{i=1}^{m}a_{ij}\theta_i)(u_1,u_2)e_j$$
and so  $\phi\mu_j-\sum_{i=1}^{m}a_{ij}\theta_i \in B_g(\mathbb U \times \mathbb U, {\bf k})$
for $j=1,\ldots,m$.

Therefore
$\left[ \phi \mu_j \right]=\sum_{i=1}^{m}a_{ij}\left[\theta_i\right]$ for every $j=1,...,m$ and
we get that the sets $\{\left[ \phi\theta_i \right]\}_{i\leq m}$ and $\{\left[ \mu_i \right]\}_{i\leq m}$ span the same subspace of
$H^2(\mathbb U, {\bf k})$.

Conversely,  suppose $\{\left[ \phi\theta_i \right]\}_{i\leq m}$
spans the same subspace of $H^2(\mathbb U, {\bf k})$ than
$\{\left[ \mu_i \right]\}_{i\leq m}$. Then there exist an invertible  matrix
$(a_{ij})_{ij}$ such that
$\left[\phi \mu_j\right]= \sum_{i=1}^{m}a_{ij}\left[\theta_i\right]$.
Define the linear isomorphism:
\begin{align*}
\psi:\mathbb W&\rightarrow \mathbb W \\
\psi(e_i)&\mapsto \sum_{j=1}^{m}a_{ij}e_j.
\end{align*}

Then
$\psi \theta(u_1,u_2)=\sum_{i=1}^{m}\sum_{j=1}^{m}a_{ij}\theta_i(u_1,u_2)e_j$ and also
$$(\phi\mu-\psi\theta)(u_1,u_2)=\sum_{j=1}^{m}(\phi\mu_j-\sum_{i=1}^{m}a_{ij}\theta_i)(u_1,u_2)e_j.$$

From here $\left[\phi \mu\right]=\left[\psi \theta\right]$ and, by Lemma
\ref{l24}, we conclude that the bilinear pairs $(\mathbb U\oplus \mathbb W, g_{\theta})$ and $(\mathbb U\oplus \mathbb W,g_{\mu})$ are equivalent.
\end{Proof}


Let us introduce some of notation.

\begin{notation}
Let $(\mathbb U,g)$ be a bilinear pair. We denote by
$$\hbox{${\mathcal{Aut}} (\mathbb U,g):=\{\phi:\mathbb U \to \mathbb U \in GL(\mathbb U)$ such that $g(\phi(u_1),\phi(u_2))=\phi(g(u_1,u_2))$ for any $u_1,u_2 \in \mathbb U\}$.}$$
Observe that ${\mathcal{Aut}} (\mathbb U,g)$ is a subgroup of $GL(\mathbb U)$ that will be called the {\rm equivalence group of} $(\mathbb U,g)$ .
\end{notation}

We recall that for a finite-dimensional linear space $\mathbb V$, the {\it Grassmannian} $Grass_s(\mathbb V)$ is the set of all $s$-dimensional linear subspaces of $\mathbb V$.

\begin{notation}
Let $(\mathbb U,g)$ be a bilinear pair. For any $1 \leq s \leq {\rm dim}(\mathbb U)^2-{\rm dim}(g(\mathbb U,\mathbb U))$ we denote by
$$T_s(\mathbb U)=\{ \langle [\theta_1],\ldots,[\theta_s]\rangle \in Grass_s(H^2(\mathbb U,{\bf k})): (\cap_{i=1}^{s}{\rm rad}(\theta_i))\cap {\rm rad}(g)=0  \}$$
where $Grass_s(H^2(\mathbb U,{\bf k}))$ is the Grassmannian of subspaces of dimension $s$ in
$H^2(\mathbb U,{\bf k})$.
\end{notation}

The next result is of easy verification.

\begin{lemma}
The set $T_s(\mathbb U)$ is stable under the  action of ${\mathcal A}(\mathbb U,g)$.
\end{lemma}

\begin{definition}
Let $(\mathbb U,g)$ be a bilinear pair. A radical extension of $(\mathbb U,g)$ is a
bilinear pair $(\mathbb U\oplus \mathbb W, g_{\theta})$, where $\mathbb W$ is a vector space and
$\theta\in Bil(\mathbb U\times \mathbb U, \mathbb W)$.
\end{definition}


By Lemma \ref{l22}, any bilinear pair $(\mathbb V, f)$ of dimension $n$ with ${\rm dim}({\rm rad}(f))=m \neq 0$ is equivalent to one the form $(\mathbb U\oplus \mathbb W, g_{\theta})$ where $(\mathbb U, g)$ is a bilinear pair of dimension $m-n$, $\mathbb W$ is a linear space with ${\rm dim}(\mathbb W)={\rm dim}({\rm rad}(f))$,  and $\theta: \mathbb U \times \mathbb U \to \mathbb W$ is
such that ${\rm rad}(g)\cap {\rm rad}(\theta)=0$. That is $(\mathbb V, f)$ is equivalent to a radical extension
 of $(\mathbb U,g)$.
Also observe that,  by Lemma \ref{l211}, there is a correspondence between the set of equivalent classes of radical extensions
(without radical component) of $(\mathbb U,g)$ and the set of orbits
of $T_s(\mathbb U)$ under the action of ${\mathcal A}(\mathbb U,g)$.

\bigskip

Summarizing,  we have obtained a procedure   to classify (up to equivalence) all of  the bilinear pairs of dimension $n$ with
${\rm dim}({\rm rad}(f))=m \neq 0$, in case we  know the classification  of the  bilinear
pairs $(\mathbb U,g)$ of dimension $n-m$. This follows  three steps for any pair $(\mathbb U,g)$:

\medskip

{\centerline{{\bf PROCEDURE 1}}}

\begin{itemize}
\item[(i)] Compute  ${\mathcal A}(\mathbb U,g)$, $H^2_g(\mathbb U, {\bf k})$  and ${{\rm rad}} (g)$.
\item[(ii)] Determine the set of ${\mathcal A}(\mathbb U,g) $-orbits on $T_s(\mathbb U)$, where $1 \leq s \leq (n-m)^2-{\rm dim}(g(\mathbb U,\mathbb U))$.
\item[(iii)] For each orbit, construct the bilinear pair  corresponding to a representative of it.
\end{itemize}

\bigskip

We will apply the above procedure to  classify the $n$-dimensional bilinear pairs with a radical of dimension $n-2$. So we begin by showing the classification of the two-dimensional bilinear pairs.

\subsection{Classification of $2$-dimensional bilinear pairs}

Thanks to \cite{kv17} we have  the classification (up to isomorphism) of two-dimensional
algebras, which is also the classification (up to equivalence) of two-dimensional bilinear
pair.

To give this  classification we need to introduce some notation.
Let us consider the action of the cyclic group $C_2=\langle \rho\mid \rho^2\rangle$ on ${\bf k}$ defined by the equality ${}^{\rho}\alpha=-\alpha$ for $\alpha\in{\bf k}$.
Let us fix some set of representatives of orbits under this action and denote it by ${\bf k_{\geq 0}}$. For example, if ${\bf k}=\mathbb{C}$, then one can take $\mathbb{C}_{\ge 0}=\{\alpha\in\mathbb{C}\mid Re(\alpha)>0\}\cup\{\alpha\in\mathbb{C}\mid Re(\alpha)=0,Im(\alpha)\ge 0\}$.

Let us also consider the action of $C_2$ on ${\bf k}^2$ defined by the equality ${}^{\rho}(\alpha,\beta)=(1-\alpha+\beta,\beta)$ for $(\alpha,\beta)\in{\bf k}^2$.
Let us fix some set of representatives of orbits under this action and denote it by $\mathcal{U}$. Let us also define $\mathcal{T}=\{(\alpha,\beta)\in{\bf k}^2\mid \alpha+\beta=1\}$.

Given $(\alpha,\beta,\gamma,\delta)\in{\bf k}^4$, we define $\mathcal{D}(\alpha,\beta,\gamma,\delta)=(\alpha+\gamma)(\beta+\delta)-1$.
We define $\mathcal{C}_1(\alpha,\beta,\gamma,\delta)=(\beta,\delta)$, $\mathcal{C}_2(\alpha,\beta,\gamma,\delta)=(\gamma,\alpha)$, and $\mathcal{C}_3(\alpha,\beta,\gamma,\delta)=\left(\frac{\beta\gamma-(\alpha-1)(\delta-1)}{\mathcal{D}(\alpha,\beta,\gamma,\delta)},\frac{\alpha\delta-(\beta-1)(\gamma-1)}{\mathcal{D}(\alpha,\beta,\gamma,\delta)}\right)$ for $(\alpha,\beta,\gamma,\delta)$ such that $\mathcal{D}(\alpha,\beta,\gamma,\delta)\not=0$. Let us consider the set $\left\{\big(\mathcal{C}_1(\Gamma),\mathcal{C}_2(\Gamma),\mathcal{C}_3(\Gamma)\big)\mid \Gamma\in{\bf k}^4, \mathcal{D}(\Gamma)\not=0,\mathcal{C}_1(\Gamma),\mathcal{C}_2(\Gamma)\not\in \mathcal{T}\right\}\subset ({\bf k}^2)^3$.
One can show that the symmetric group $S_3$ acts on this set by the equality ${}^{\sigma}\big(\mathcal{C}_1(\Gamma),\mathcal{C}_2(\Gamma),\mathcal{C}_3(\Gamma)\big)=\big(\mathcal{C}_{\sigma^{-1}(1)}(\Gamma),\mathcal{C}_{\sigma^{-1}(2)}(\Gamma),\mathcal{C}_{\sigma^{-1}(3)}(\Gamma)\big)$ for $\sigma\in S_3$. Note that there exists a set of representatives of orbits under this action $\mathcal{\tilde V}$ such that if $(\mathcal{C}_1,\mathcal{C}_2,\mathcal{C}_3)\in \mathcal{\tilde V}$ and $\mathcal{C}_1\not=\mathcal{C}_2$, then $\mathcal{C}_3\not=\mathcal{C}_1,\mathcal{C}_2$. Let us fix such $\mathcal{\tilde V}$ and define
$$
\mathcal{V}=\{\Gamma\in{\bf k}^4\mid \mathcal{D}(\Gamma)\not=0; \mathcal{C}_1(\Gamma),\mathcal{C}_2(\Gamma)\not\in \mathcal{T}, \big(\mathcal{C}_1(\Gamma),\mathcal{C}_2(\Gamma),\mathcal{C}_3(\Gamma)\big)\in\mathcal{\tilde V}\}.
$$
For $\Gamma\in\mathcal{V}$, we also define $\mathcal{C}(\Gamma)=\{\mathcal{C}_1(\Gamma),\mathcal{C}_2(\Gamma),\mathcal{C}_3(\Gamma)\}\subset{\bf k}^2$.

Let us consider the action of the cyclic group $C_2$ on ${\bf k}^*\setminus \{1\}$ defined by the equality ${}^{\rho}\alpha=\alpha^{-1}$ for $\alpha\in{\bf k}^*\setminus \{1\}$.
Let us fix some set of representatives of orbits under this action and denote it by ${\bf k_{>1}^*}$. For example, if ${\bf k}=\mathbb{C}$, then one can take $\mathbb{C}_{>1}^*=\{\alpha\in\mathbb{C}^*\mid |\alpha|>1\}\cup\{\alpha\in\mathbb{C}^*\mid |\alpha|=1,0<arg(\alpha)\le \pi\}$. For $(\alpha,\beta,\gamma)\in {\bf k}^2\times{\bf k}^*_{>1}$ we define $$\mathcal{C}(\alpha,\beta,\gamma)=\left\{\big(\alpha\gamma,(1-\alpha)\gamma\big),\left(\frac{\beta}{\gamma},\frac{1-\beta}{\gamma}\right)\right\}\subset{\bf k}^2.$$

In \hyperref[tab4]{Table 4}, we summarize the classification up to equivalence  of two-dimensional bilinear pairs $(\mathbb{U},g)$ given in \cite{kv17}. Fix a basis $\{e_1,e_2\}$ of $\mathbb{U}$, then the first column will denote the given name to  bilinear pairs, the second column the description of  bilinear maps $g(e_i,e_j):=e_ie_j$, $i, j \in \{1,2\}$ and the third column is devoted to  show  the corresponding equivalence groups,  which can be also found in  \cite{kv17}.

\begin{remark}\label{remrad}
Let $({\mathbb U},g)$ be a $2$-dimensional bilinear pair  with a nonzero radical,
then $({\mathbb U},g)$ is equivalent to  ${\bf A}_3 $, ${\bf D}_2(0,0)$ or $\mathfrak{N}_2$.
\end{remark}

\section{Proof of the Main theorem}

In this section we will apply Procedure 1 to prove our classification theorem of $n$-dimensional bilinear pairs with a radical of dimension $n-2$ (Theorem \ref{maintheo}). Hence, we will consider each two-dimensional bilinear pair $(\mathbb{U},g)$ in \hyperref[tab4]{Table 4}. For each one we will follow the three steps described in Procedure 1.

\medskip

First, we recall that the set $\{\Delta_{ij} \}_{1\leq i,j\leq 2}$, where each
$\Delta_{ij}:\mathbb U\times \mathbb U \rightarrow {\bf k}$ is defined by
$\Delta_{ij}(e_i, e_j)=1$ and $\Delta_{ij}=0$ otherwise, is a basis of
$Bil(\mathbb  U\times \mathbb U,{\bf k})$. Also, for any bilinear pair $(\mathbb U,g)$, we will denote by $({\mathbb U}^{op}, g^{op})$ the {\it opposite bilinear pair} given by  ${\mathbb U}^{op}={\mathbb U}$
and $g^{op}(x,y):= g(y,x)$ for any $x,y \in {\mathbb U}^{op}$. In this context, we have  the  following remark which becomes useful when applying the procedure for some cases.

\begin{remark}\label{Wish}
Let $(\mathbb U,g)$ be a $k$-dimensional bilinear pair  and denote by $\{ ({\mathbb U}_{n,i}, g_{n,i})\}_{i \in I}$  the set of equivalence classes of
$n$-dimensional radical extensions of $(\mathbb U,g).$
Then there is a one-to-one correspondence between $\{ ({\mathbb U}_{n,i}, g_{n,i})\}_{i \in I}$ and $\{ (({\mathbb U}^{op})_{n,j}, (g^{op})_{n,j})\}_{j \in J}$.
\end{remark}
\medskip

\subsection{Bilinear pairs ${\bf A}_1 (\alpha)$, $\alpha\in{\bf k}$}

We begin by noting that the bilinear pair  $({\bf A}_1 (0),g)$ is equivalent to the bilinear pair  $({\bf A}_1 (1)^{op}, g^{op})$. Hence (see Remark \ref{Wish}), we will now consider bilinear pairs $({\bf A}_1 (\alpha),g)$ with $\alpha \neq 0.$

\medskip

We start by calculating $H^2({\bf A}_1 (\alpha),{\bf k})$ and $T_s({\bf A}_1 (\alpha))$.

\medskip

It is easy to see that
$\delta e_1^*= \Delta_{11}$ and that  $\delta  e_2^*=\Delta_{11}+\alpha \Delta_{12} + (1-\alpha) \Delta_{21}.$ From here
 Lemma \ref{l28} gives us  $$B_g({\bf A}_1 (\alpha)\times {\bf A}_1 (\alpha), {\bf k})=\langle \delta e_1^*, \delta e_2^* \rangle= \langle \Delta_{11}, \Delta_{12} + \frac{1-\alpha}{\alpha} \Delta_{21}\rangle.$$ Then we have that  ${\rm dim} (H^2({\bf A}_1 (\alpha),{\bf k}))=4-2=2$ and that $$H^2({\bf A}_1 (\alpha),{\bf k})=\langle  \Delta_{ij}:i,j \in \{1,2\} \rangle / \langle \Delta_{11}, \Delta_{12} + \frac{1-\alpha}{\alpha} \Delta_{21}\rangle =
             \langle [ \Delta_{21}], [\Delta_{22}]\rangle.  $$

\medskip
Since ${\rm rad}(g)=0$ then
 $$
\hbox{$T_s({\bf A}_1 (\alpha))=Grass_s(H^2({\bf A}_1 (\alpha),{\bf k}))$ for $s \in \{1,2\}$.}$$

\medskip

Now, we calculate the orbits of each space in $T_s({\bf A}_1 (\alpha))$ under the
action of the group ${\mathcal A}({\bf A}_1 (\alpha),g)$.

The action of an element $ \left(
                             \begin{array}{cc}
                               1 & 0 \\
                               x & 1 \\
                             \end{array}
                           \right) \in
{\mathcal A}({\bf A}_1 (\alpha),g)$ on a subspace
$\langle a [\Delta_{21}] +b [\Delta_{22}]   \rangle  \in T_1({\bf A}_1 (\alpha))$
is the following:
$$\langle (bx^2+ax) [\Delta_{11}] + bx [\Delta_{12}] + (bx+a)[\Delta_{21}] + b[\Delta_{22}] \rangle =
\langle \Big(a+\Big(1-\dfrac{1-\alpha}{\alpha}\Big)bx\Big) [\Delta_{21}] + b [\Delta_{22}] \rangle.$$

Since we can write $T_1({\bf A}_1 (\alpha\neq 0))=
\left\{ \langle [\Delta_{21}]   \rangle\right\} \bigcup \left\{ \langle \lambda [\Delta_{21}] + [\Delta_{22}]  \rangle: \lambda \in {\bf k}\right\} $ we are going to distinguish two cases.

\begin{enumerate}
    \item[$\bullet$] $b=0$, the spaces in the orbit of
$\langle  a [\Delta_{21}]  \rangle$
are
$\langle  a [\Delta_{21}]  \rangle  =
\langle  [\Delta_{21}]  \rangle,$ for $x\in {\bf k}$.

\item[$\bullet$]  $b\not=0$, we denote $\lambda=a/b.$
The orbit of
$\langle \lambda [\Delta_{21}] + [\Delta_{22}] \rangle$
is
$$\left\{\langle\Big(\lambda +\Big(1-\dfrac{1-\alpha}{\alpha}\Big)x\Big) [\Delta_{21}] + [\Delta_{22}] \rangle,
\mbox{ for }x\in {\bf k}\right\}.$$
\begin{enumerate}

\item[$\diamond$] if $\alpha\neq \frac{1}{2},$ we can choose a good representation:
for $x=-\lambda \Big(1-\dfrac{1-\alpha}{\alpha}\Big)^{-1}$ we have
$\langle  [\Delta_{22}]   \rangle.$

\item[$\diamond$] if $\alpha=\frac{1}{2},$ we have one stable orbit for each element $\lambda \in {\bf k}$ we have
a representation
$\langle \lambda [\Delta_{21}] + [\Delta_{22}]  \rangle.$

\end{enumerate}

\end{enumerate}

\medskip

 Since $T_2({\bf A}_1 (\alpha))=H^2({\bf A}_1 (\alpha), {\bf k})$ we just have one orbit in the action which is the whole $H^2({\bf A}_1 (\alpha),{\bf k}).$

At last, we compute each class of equivalence of radical extensions of $({\bf A}_1 (\alpha), g)$ using its
correspondence with the set of  orbits in $T_s({\bf A}_1 (\alpha))$, $s \in \{1,2\}$,  under the action of the group ${\mathcal A}({\bf A}_1 (\alpha),g)$. We will also take into account the above mentioned equivalence between the bilinear pairs
$({\bf A}_1 (0),g)$ and  $({\bf A}_1 (1)^{op}, g^{op})$.

$$\begin{array}{lll lll}
({\bf A}_1 (\alpha\neq0))_{3,1} &:&
e_1 e_1 = e_1+e_2 & e_1 e_2=\alpha e_2 & e_2 e_1=(1-\alpha)e_2+ e_3 & e_2 e_2=0 \\
({\bf A}_1 (0))_{3,1} &:&
e_1 e_1 = e_1+e_2 & e_1 e_2=e_3 & e_2 e_1=e_2 & e_2 e_2=0 \\
({\bf A}_1 (\alpha\neq\frac{1}{2}))_{3,2} &:&
e_1 e_1 = e_1+e_2 & e_1 e_2=\alpha e_2 & e_2 e_1=(1-\alpha)e_2 & e_2 e_2=e_3 \\
({\bf A}_1 (\frac{1}{2}))_{3,2}(\lambda\in {\bf k}) &:&
e_1 e_1 = e_1+e_2 & e_1 e_2=\frac{1}{2} e_2 & e_2 e_1=\frac{1}{2} e_2 +\lambda e_3 & e_2 e_2=e_3 \\
({\bf A}_1 (\alpha\neq0))_{4,3} &:&
e_1 e_1 = e_1+e_2 & e_1 e_2=\alpha e_2 & e_2 e_1=(1-\alpha)e_2+e_3 & e_2 e_2=e_4 \\
({\bf A}_1 (0))_{4,3} &:&
e_1 e_1 = e_1+e_2 & e_1 e_2= e_3 & e_2 e_1=e_2 & e_2 e_2=e_4 \\

 \end{array}$$

\subsection{Bilinear pair ${\bf A}_2 $}

Again, we calculate $H^2({\bf A}_2 ,{\bf k})$ and $T_s({\bf A}_2 )$.
It is easy to see that $\delta  e_1^*=0$ and that
$$\delta  e_2^*=\Delta_{11}+\Delta_{12}-\Delta_{21}.$$ Then
 by Lemma \ref{l28} we have that ${\rm dim}(H^2({\bf A}_2,{\bf k}))=3$ and that

$$H^2({\bf A}_2 ,{\bf k})=\langle  \Delta_{ij}:i,j \in \{1,2\} \rangle / \langle \Delta_{11}+ \Delta_{12} - \Delta_{21}\rangle =
             \langle              [\Delta_{11}], [\Delta_{21}], [\Delta_{22}] \rangle.$$

             \medskip

Since ${\rm rad}(g)=0$ then $$\hbox{$T_s({\bf A}_2 )=Grass_s(H^2({\bf A}_2 ,{\bf k}))$ for $s \in \{1,2,3\}$.}$$

Additionally, we calculate the orbits of each subspace of $T_s({\bf A}_2 )$ under
the action of the group ${\mathcal A}({\bf A}_2,g )$.

\begin{enumerate}
\item[$\star$] $T_1({\bf A}_2 )=
\left\{\langle [ \Delta_{11} ]  \rangle\right\} \cup
\left\{\lambda\langle [ \Delta_{11} ]  + [ \Delta_{21} ]  \rangle : \lambda \in {\bf k}\right\} \cup
\left\{\langle  \lambda [\Delta_{11}] +  \mu [\Delta_{21}] + [\Delta_{22}]  \rangle : \lambda, \mu \in {\bf k}   \right\}.$

The action of an element $\left(
                             \begin{array}{cc}
                               1 & 0 \\
                               x & 1 \\
                             \end{array}
                           \right) \in {\mathcal A}({\bf A}_2,g )$ on a subspace
$\langle
 a[\Delta_{11}] + b[\Delta_{21}] + c[\Delta_{22}] \rangle  \in T_1({\bf A}_2 )$
is the following:
$$\langle
(a+bx+cx^2)[\Delta_{11}] + cx [\Delta_{12}]+(b+cx)[\Delta_{21}]+ c[\Delta_{22}] \rangle =
\langle (a+(b-c)x+cx^2)[\Delta_{11}]+ (b+2cx) [\Delta_{21}] +  c[\Delta_{22}]\rangle.$$

From here, we can consider  $3$ situations:
\begin{enumerate}
\item[$\bullet$] $b=0$ and $c=0.$
The orbit of the space
$\langle a [\Delta_{11}] \rangle$ is trivial.

\item[$\bullet$] $b\not=0$ and $c=0$, we denote $\lambda=a/b.$
The orbit space of
$\langle  \lambda [\Delta_{11}] + [\Delta_{21}] \rangle $
is the following
$\langle (\lambda+x)[\Delta_{11}]+ [\Delta_{21}]\rangle.$
For $x=-\lambda$ we have the following representative $\langle [\Delta_{21}]\rangle.$

\item[$\bullet$] $c\not=0$, we denote $\lambda=a/c$ and $\mu= b/c.$
The orbit space of
$\langle  \lambda [\Delta_{11}] + \mu [\Delta_{21}] + [\Delta_{22}] \rangle$ is
$$\{ \langle(\lambda +(\mu-1)x+x^2)[\Delta_{11}] + (\mu+2x)[\Delta_{21}] + [\Delta_{22}]   \rangle: x\in {\bf k} \}.$$
For $x=-\frac{\mu}{2},$
we have the following representative $\langle \lambda' [\Delta_{11}] +[\Delta_{22}] \rangle$, where $\lambda'=\lambda -\frac{(\mu-1)\mu}{2}+\frac{\mu^2}{4}$. Therefore,  we have one stable orbit for each $\lambda' \in {\bf k}$.
\end{enumerate}

    \item[$\star$]
$T_2({\bf A}_2)= \left\{ \langle  [\Delta_{21}], [\Delta_{22}]  \rangle   \right\} \cup \left\{ \langle  [\Delta_{11}]+\lambda[\Delta_{22}], [\Delta_{21}]  \rangle: \lambda \in {\bf k}   \right\} \cup
\left\{ \langle  [\Delta_{11}]+\lambda[\Delta_{21}] , [\Delta_{22}]  \rangle: \lambda \in {\bf k}   \right\} \cup
\left\{ \langle  \lambda[\Delta_{21}]+[\Delta_{22}], [\Delta_{11}]  \rangle: \lambda \in {\bf k}^*   \right\} \cup
\left\{ \langle[\Delta_{11}]+\lambda[\Delta_{21}],  \mu [\Delta_{21}]+[\Delta_{22}]   \rangle: \lambda, \mu \in {\bf k}^*. \right\}$

\noindent The action of ${\mathcal A}({\bf A}_2,g )$ on a subspace $\mathfrak V$ spanned by
$\alpha_1[\Delta_{11}]+\alpha_2[\Delta_{21}]+\alpha_3[\Delta_{22}]$ and
$\beta_1 [\Delta_{11}]+\beta_2 [\Delta_{21}]+\beta_3 [\Delta_{22}]$ is the following:
$$\langle (\alpha_1+(\alpha_2-\alpha_3)x+\alpha_3  x^2)[\Delta_{11}]+
(\alpha_2 +2\alpha_3 x) [\Delta_{21}] +  \alpha_3 [\Delta_{22}],
(\beta_1+(\beta_2-\beta_3)x+\beta_3  x^2)[\Delta_{11}]+
(\beta_2 +2\beta_3 x) [\Delta_{21}] +  \beta_3 [\Delta_{22}] \rangle.$$

From here, there are $5$ cases to distinguish:
\begin{enumerate}
    \item[$\bullet$] The orbit of $\langle  [\Delta_{21}], [\Delta_{22}]  \rangle$ is $\left\{ \langle[\Delta_{11}]+\frac{1}{x}[\Delta_{21}],  (1+x) [\Delta_{21}]+[\Delta_{22}]   \rangle: x\in {\bf k} \right\}\cup \left\{\langle  [\Delta_{21}], [\Delta_{22}]  \rangle\right\}$.

    \item[$\bullet$] The orbit of $ \langle  [\Delta_{11}]+\lambda[\Delta_{22}], [\Delta_{21}]  \rangle$ for $\lambda\neq0$ is $\left\{ \langle[\Delta_{11}]+\frac{1}{x}[\Delta_{21}],  (x+1-\frac{1}{\lambda x}) [\Delta_{21}]+[\Delta_{22}]   \rangle: x\in {\bf k} \right\}\cup \left\{\langle  [\Delta_{11}]+\lambda[\Delta_{22}], [\Delta_{21}]  \rangle\right\}$. It is easy to prove that  we have one stable orbit for each $\lambda \in {\bf k}^*$.

    For $\lambda=0$, we have the trivial orbit: $\langle[\Delta_{11}], [\Delta_{21}]  \rangle$ .

    \item[$\bullet$] The subspaces $\langle  [\Delta_{11}]+\lambda'[\Delta_{21}] , [\Delta_{22}]  \rangle$ for $\lambda'\neq0$ belong to the previous case if $\lambda=\frac{\lambda'}{1+\frac{1}{\lambda'}}$, chosing $x=1/\lambda'$.

    If $\lambda'=0$ then the orbit is $\left\{ \langle[\Delta_{11}], 2x [\Delta_{21}]+[\Delta_{22}] \rangle: x\in {\bf k}\right\} $ and we can choose a good representation: for $x=0$ we have $\langle  [\Delta_{11}], [\Delta_{22}]  \rangle$.

    \item[$\bullet$] The subespaces $ \langle  \lambda[\Delta_{21}]+[\Delta_{22}], [\Delta_{11}]  \rangle$ for $\lambda\neq0$ belong to the orbit of $\langle  [\Delta_{11}], [\Delta_{22}]  \rangle$.

    \item[$\bullet$] Also, the subespaces $\langle[\Delta_{11}]+\lambda[\Delta_{21}],  \mu [\Delta_{21}]+[\Delta_{22}]   \rangle$ belong to the cases one and two.

\end{enumerate}

    \item[$\star$]  $T_3({\bf A}_2 )=H^2({\bf A}_2 , {\bf k}).$

\end{enumerate}

Finally, we can present each radical extension, up to equivalence,  of $({\bf A}_2 , g)$ using its
correspondence with an orbit in $T_s({\bf A}_2 )$, $s \in \{1,2,3\}$:

$$\begin{array}{lll lll}

({\bf A}_2 )_{3,1} &:&    e_1 e_1 = e_2+e_3 & e_1 e_2= e_2 & e_2 e_1=-e_2  & e_2 e_2=0\\

({\bf A}_2 )_{3,2} &:&    e_1 e_1 = e_2 & e_1 e_2= e_2 & e_2 e_1=-e_2+e_3  & e_2 e_2=0\\

({\bf A}_2 )_{3,3}(\lambda\in {\bf k}) &:&    e_1 e_1 = e_2+\lambda e_3 & e_1 e_2= e_2 & e_2 e_1=-e_2  & e_2 e_2=e_3 \\

({\bf A}_2 )_{4,4} &:&    e_1 e_1 = e_2 & e_1 e_2= e_2 & e_2 e_1=-e_2+e_3  & e_2 e_2=e_4\\

({\bf A}_2 )_{4,5}(\lambda\in {\bf k}) &:&    e_1 e_1 = e_2+e_3 & e_1 e_2= e_2 & e_2 e_1=-e_2+e_4  & e_2 e_2=\lambda e_3\\

({\bf A}_2 )_{4,6} &:&    e_1 e_1 = e_2+e_3 & e_1 e_2= e_2 & e_2 e_1=-e_2 & e_2 e_2=e_4\\

({\bf A}_2 )_{5,7}  &:&    e_1 e_1 = e_2+ e_3 & e_1 e_2= e_2 & e_2 e_1=-e_2+e_4  & e_2 e_2=e_5

 \end{array}$$

\subsection{Bilinear pair ${\bf A}_3 $}
At first, we get $H^2({\bf A}_3 ,{\bf k})$  and $T_s({\bf A}_3 )$.
It is easy to see that
$$\hbox{$\delta  e_1^*=0$ and that $\delta  e_2^*=\Delta_{11}$ }.$$
Then  by Lemma \ref{l28},
$$Basis(H^2({\bf A}_3 ,{\bf k}))=
             \left\{ [\Delta_{12}],  [\Delta_{21}], [\Delta_{22}] \right\}.$$

Since ${\rm rad}(g)=\langle e_2 \rangle$ and $e_2 \notin {\rm rad}(\theta)$ for $[\theta] \in H^2({\bf A}_3,{\bf k})$ and  $[\theta]\neq0$ then $T_s({\bf A}_3) = Grass_s(H^2({\bf A}_3 , {\bf k})).$

Now, the action of ${\mathcal A}({\bf A}_3,g )$ on a subspace $\langle a[\Delta_{12}] +b [\Delta_{21}]+ c[\Delta_{22}] \rangle$ is the following:
$$\langle (ax^3+cx^2y)[\Delta_{12}]+ (bx^3+cx^2y)[\Delta_{21}]+ cx^4 [\Delta_{22}] \rangle.$$

We denote the orbit of a subspace $\mathbb V$ by $O(\mathbb V)$. From here we omit the calculus and we just present each $T_s$ as a joint of their orbits.
\begin{enumerate}
    \item[$\star$] $T_1({\bf A}_3 )=
    O(\langle [\Delta_{22}] \rangle) \cup
    O(\langle [\Delta_{21}]+[\Delta_{22}]\rangle) \cup
    O(\langle [\Delta_{21}]\rangle) \cup O(\{\langle [\Delta_{12}]+\lambda[\Delta_{21}]\rangle:\lambda\in{\bf k}\}).$

    \item[$\star$] $T_2({\bf A}_3 )=$ $ O( \langle  [\Delta_{22}], [\Delta_{21}]  \rangle)  \cup   O( \langle  [\Delta_{21}], [\Delta_{12}] \rangle) \cup O(\{ \langle [\Delta_{22}], [\Delta_{12}]+\lambda[\Delta_{21}] \rangle:\lambda\in{\bf k}\})).$

    \item[$\star$]  $T_3({\bf A}_3 )=H^2({\bf A}_3 , {\bf k})$.

\end{enumerate}

At last, we compute each radical extension of $({\bf A}_3 , g)$ using its
correspondence with an orbit in $T_s({\bf A}_3 )$.
$$\begin{array}{ll llll}
({\bf A}_3 )_{3,1} &:& e_1 e_1 = e_2 & e_1 e_2=0 & e_2 e_1=0 & e_2 e_2=e_3 \\
({\bf A}_3 )_{3,2} &:& e_1 e_1 = e_2 & e_1 e_2=0 & e_2 e_1= e_3 & e_2 e_2=e_3 \\
({\bf A}_3 )_{3,3} &:& e_1 e_1 = e_2 & e_1 e_2=0 & e_2 e_1=e_3 & e_2 e_2=0 \\
({\bf A}_3 )_{3,4}(\lambda\in {\bf k}) &:& e_1 e_1 = e_2 & e_1 e_2=e_3 & e_2 e_1=\lambda e_3 & e_2 e_2=0 \\
({\bf A}_3 )_{4,5} &:& e_1 e_1 = e_2 & e_1 e_2=0 & e_2 e_1=e_4 & e_2 e_2=e_3 \\
({\bf A}_3 )_{4,6} &:& e_1 e_1 = e_2 & e_1 e_2=e_4 & e_2 e_1=e_3 & e_2 e_2=0 \\
({\bf A}_3 )_{4,7} (\lambda\in {\bf k}) &:& e_1 e_1 = e_2 & e_1 e_2=e_4 & e_2 e_1=\lambda e_4 & e_2 e_2=e_3 \\
({\bf A}_3 )_{5,8} &:& e_1 e_1 = e_2 & e_1 e_2=e_3 & e_2 e_1=e_4 & e_2 e_2=e_5

\end{array}$$

\subsection{Bilinear pairs ${\bf A}_4 (\alpha)$, $\alpha\in{\bf k_{\geq 0}}$}

To start with, we calculate $H^2({\bf A}_4 (\alpha),{\bf k})$  and $T_s({\bf A}_4 (\alpha))$.
It is easy to see that
$$\delta  e_1^*= \alpha \Delta_{11} +\Delta_{12} - \Delta_{21}, \ \delta e_2^*=\Delta_{11}+\alpha \Delta_{12},$$
and by Lemma \ref{l28},
$$Basis(H^2({\bf A}_4 (\alpha), {\bf k})) = \{ [\Delta_{12}], [\Delta_{22}] \}.$$
Since ${\rm rad}(g)=0,$ then $T_s({\bf A}_4 (\alpha)) = Grass_s(H^2({\bf A}_4 (\alpha), {\bf k})).$

Moreover, we calculate the orbits of each space in $T_s({\bf A}_4 (\alpha))$ under the action of ${\mathcal A}({\bf A}_4 (\alpha),g)$. Note that the action of ${\mathcal A}({\bf A}_4 (\alpha),g)$ on a subspace
$\langle  a[\Delta_{12}] +  b[\Delta_{22}] \rangle $
is $\langle  \gamma [\Delta_{12}] +  b[\Delta_{22}] \rangle$
($\gamma=a$ for $\alpha \neq 0$ and $\gamma =\pm a$ for $\alpha=0$).

It is easy to see that

\begin{enumerate}
    \item[$\star$] $T_1({\bf A}_4 (\alpha\neq0))=
   O( \langle [\Delta_{22}] \rangle) \cup
    O(\{ \langle [\Delta_{12}]+\lambda [\Delta_{22}]\rangle: \lambda \in {\bf k} \}).$

    \item[$\star$] $T_1({\bf A}_4 (0))=
    O(\langle [\Delta_{22}] \rangle) \cup
    O(\{ \langle [\Delta_{12}]+\lambda [\Delta_{22}]\rangle: \lambda \in {\bf k}_{\geq 0}\}).$

    \item[$\star$]  $T_2({\bf A}_4 (\alpha))=H^2({\bf A}_4 (\alpha), {\bf k}).$

\end{enumerate}

To sum up, we present the radical extensions of ${\bf A}_4(\alpha)$, up to equivalence.
$$\begin{array}{ll llll}
({\bf A}_4 (\alpha))_{3,1} &:& e_1 e_1 = \alpha e_1+e_2 & e_1 e_2= e_1+\alpha e_2 & e_2 e_1=-e_1  & e_2 e_2= e_3 \\
({\bf A}_4 (\alpha\neq 0))_{3,2}(\lambda\in{\bf k}) &:& e_1 e_1 = \alpha e_1+e_2 & e_1 e_2= e_1+\alpha e_2 + e_3 & e_2 e_1=-e_1   & e_2 e_2=\lambda e_3 \\
({\bf A}_4 (0))_{3,2}(\lambda\in{\bf k}_{\geq 0}\}) &:& e_1 e_1 = e_2 & e_1 e_2= e_1 + e_3 & e_2 e_1=-e_1   & e_2 e_2=\lambda e_3 \\
({\bf A}_4 (\alpha))_{4,3}&:& e_1 e_1 = \alpha e_1+e_2 & e_1 e_2= e_1+\alpha e_2+ e_3 & e_2 e_1=-e_1  & e_2 e_2= e_4
\end{array}$$

\subsection{Bilinear pairs ${\bf B}_1(\alpha)$, $\alpha\in{\bf k}$}
Firstly, we calculate $H^2({\bf B}_1(\alpha),{\bf k})$  and $T_s({\bf B}_1(\alpha))$.
It is easy to see that
$$\delta  e_1^*=(1- \alpha) \Delta_{12} + \alpha \Delta_{21}, \ \delta e_2^*=\Delta_{12}- \Delta_{21},$$
and by Lemma \ref{l28},
$$Basis(H^2({\bf B}_1(\alpha),{\bf k}))=
             \left\{ [\Delta_{11}], [\Delta_{22}] \right\}.$$
Since ${\rm rad}(g)=0$ then $T_s({\bf B}_1(\alpha))=Grass_s(H^2({\bf B}_1(\alpha),{\bf k}))$.

Secondly, we calculate the orbits of each space in $T_s({\bf B}_1(\alpha))$ under the action of ${\mathcal A}({\bf B}_1(\alpha),g)$.
This action is the identity.

It is easy to see that

\begin{enumerate}
    \item[$\star$] $T_1({\bf B}_1(\alpha))= \{ \langle[\Delta_{22}]\rangle \} \cup \{ \langle [\Delta_{11}] + \lambda [\Delta_{22}] \rangle: \lambda \in {\bf k}\}. $

    \item[$\star$] $T_2({\bf B}_1(\alpha))=H^2({\bf B}_1(\alpha), {\bf k}).$

\end{enumerate}

Finally, we compute each radical extension of $({\bf B}_1(\alpha), g)$ using its
correspondence with an orbit in $T_s({\bf B}_1(\alpha))$.
$$\begin{array}{ll llll}
({\bf B}_1(\alpha))_{3,1} &:& e_1 e_1 = 0 & e_1 e_2=(1- \alpha) e_1 + e_2  & e_2 e_1=\alpha e_1- e_2 &  e_2 e_2= e_3  \\
({\bf B}_1(\alpha))_{3,2}(\lambda\in {\bf k}) &:& e_1 e_1 = e_3 & e_1 e_2=(1- \alpha) e_1 + e_2  & e_2 e_1=\alpha e_1- e_2 &  e_2 e_2= \lambda  e_3  \\
({\bf B}_1(\alpha))_{4,3} &:& e_1 e_1 = e_3 & e_1 e_2=(1- \alpha) e_1 + e_2  & e_2 e_1=\alpha e_1- e_2 &  e_2 e_2= e_4
\end{array}
$$

\subsection{Bilinear pairs ${\bf B}_2(\alpha)$, $\alpha\in{\bf k}$}
Note that the bilinear pairs  ${\bf B}_2(0)$ and   ${\bf B}_2(1)^{op}$ are equivalent. Hence, now  we will consider only $\alpha \neq 0.$

 To begin with, we calculate $H^2({\bf B}_2(\alpha),{\bf k})$ and $T_s({\bf B}_2(\alpha))$.
It is easy to see that
$$\delta  e_1^*=(1- \alpha) \Delta_{12} + \alpha \Delta_{21}, \delta  e_2^*=0$$
and by Lemma \ref{l28},
$$Basis(H^2({\bf B}_2(\alpha),{\bf k}))= \left\{        [\Delta_{11}], [\Delta_{12}], [\Delta_{22}] \right\}.$$

Since ${\rm rad}(g)=0$ then $T_s({\bf B}_2(\alpha))=Grass_s(H^2({\bf B}_2(\alpha),{\bf k}))$.

Additionally, we calculate the orbits of each space in $T_s({\bf B}_2(\alpha))$ under the action of ${\mathcal A}({\bf B}_2(\alpha),g)$.

The orbit of a subspace
$\langle  a[\Delta_{11}]  +b[\Delta_{12}]+ c [\Delta_{22}] \rangle $by the action of ${\mathcal A}({\bf B}_2(\alpha),g)$ is the following set:
$$\{ \langle  ax^2[\Delta_{11}]  +bx[\Delta_{12}]+ c [\Delta_{22}] \rangle :x\in{\bf k}^*\}.$$

By the same method applied in the previous cases we get:

\begin{enumerate}
    \item[$\star$]
     $T_1({\bf B}_2(\alpha\neq0))=$
     $ O(\langle [\Delta_{22}] \rangle ) \cup
                        O( \langle [\Delta_{12}]+[\Delta_{22}] \rangle ) \cup
                        O( \langle [\Delta_{12}] \rangle ) \cup
                        O( \langle [\Delta_{11}]\rangle)
                        \cup
                        O( \langle [\Delta_{11}] + [\Delta_{22}] \rangle )\cup$

                        $\cup
                        O(\{ \langle [\Delta_{11}] + [\Delta_{12}]+ \lambda  [\Delta_{22}] \rangle : \lambda \in{\bf k}\}) .$
\item[$\star$]
$T_2({\bf B}_2(\alpha\neq0))=$
$ O(\langle [\Delta_{12}], [\Delta_{22}] \rangle ) \cup
O(\langle [\Delta_{11}]+[\Delta_{22}], [\Delta_{12}] \rangle ) \cup
O(\langle [\Delta_{11}], [\Delta_{12}] \rangle ) \cup
O(\langle [\Delta_{11}]+[\Delta_{12}], [\Delta_{22}] \rangle ) \cup$

$\cup O(\langle [\Delta_{11}], [\Delta_{22}] \rangle )
      \cup                         O(\{ \langle [\Delta_{11}] + \lambda  [\Delta_{12}] , [\Delta_{12}]+[\Delta_{22}] \rangle : \lambda\in{\bf k}\}) .$
\item[$\star$] $T_3({\bf B}_2(\alpha\neq0))=H^2({\bf B}_2(\alpha), {\bf k})$.
\end{enumerate}

Finally, we compute each radical extension of $({\bf B}_2(\alpha), g)$ using its
correspondence with an orbit in $T_s({\bf B}_2(\alpha))$, $s \in \{1,2,3\}$.
$$\begin{array}{ll llll}
({\bf B}_2(\alpha))_{3,1} &:&
e_1 e_1 = 0  & e_1 e_2=(1- \alpha) e_1 & e_2 e_1=\alpha e_1 &  e_2 e_2=e_3 \\
({\bf B}_2(\alpha\neq0))_{3,2} &:&
e_1 e_1 = 0  & e_1 e_2=(1- \alpha) e_1+e_3 & e_2 e_1=\alpha e_1 &  e_2 e_2=e_3 \\
({\bf B}_2(0))_{3,2} &:&
e_1 e_1 = 0  & e_1 e_2=e_1 & e_2 e_1=e_3 &  e_2 e_2=e_3 \\
({\bf B}_2(\alpha\neq0))_{3,3} &:&
e_1 e_1 = 0  & e_1 e_2=(1- \alpha) e_1+e_3 & e_2 e_1=\alpha e_1 &  e_2 e_2=0 \\
({\bf B}_2(0))_{3,3} &:&
e_1 e_1 = 0  & e_1 e_2= e_1 & e_2 e_1=e_3 &  e_2 e_2=0 \\
({\bf B}_2(\alpha))_{3,4} &:&
e_1 e_1 = e_3  & e_1 e_2=(1- \alpha) e_1 & e_2 e_1=\alpha e_1 &  e_2 e_2=0 \\
({\bf B}_2(\alpha))_{3,5} &:&
e_1 e_1 = e_3  & e_1 e_2=(1- \alpha) e_1 & e_2 e_1=\alpha e_1 &  e_2 e_2=e_3 \\
({\bf B}_2(\alpha\neq0))_{3,6}(\lambda\in {\bf k}) &:&
e_1 e_1 = e_3  & e_1 e_2=(1- \alpha) e_1+e_3 & e_2 e_1=\alpha e_1 &  e_2 e_2=\lambda e_3 \\
({\bf B}_2(0))_{3,6}(\lambda\in {\bf k}) &:&
e_1 e_1 = e_3  & e_1 e_2=e_1 & e_2 e_1=e_3 &  e_2 e_2=\lambda e_3 \\
({\bf B}_2(\alpha\neq0))_{4,7} &:&
e_1 e_1 = 0  & e_1 e_2=(1- \alpha) e_1+e_3 & e_2 e_1=\alpha e_1 &  e_2 e_2=e_4 \\
({\bf B}_2(0))_{4,7} &:&
e_1 e_1 = 0  & e_1 e_2= e_1 & e_2 e_1=e_3 &  e_2 e_2=e_4 \\
({\bf B}_2(\alpha\neq0))_{4,8} &:&
e_1 e_1 = e_3  & e_1 e_2=(1- \alpha) e_1+e_4 & e_2 e_1=\alpha e_1 &  e_2 e_2=e_3 \\
({\bf B}_2(0))_{4,8} &:&
e_1 e_1 = e_3  & e_1 e_2=e_1 & e_2 e_1=e_4 &  e_2 e_2=e_3 \\
({\bf B}_2(\alpha\neq0))_{4,9} &:&
e_1 e_1 = e_3  & e_1 e_2=(1- \alpha) e_1+e_4 & e_2 e_1=\alpha e_1 &  e_2 e_2=0 \\
({\bf B}_2(0))_{4,9} &:&
e_1 e_1 = e_3  & e_1 e_2= e_1 & e_2 e_1=e_4 &  e_2 e_2=0 \\
({\bf B}_2(\alpha\neq0))_{4,10} &:&
e_1 e_1 = e_3  & e_1 e_2=(1- \alpha) e_1+e_3 & e_2 e_1=\alpha e_1 &  e_2 e_2=e_4 \\
({\bf B}_2(0))_{4,10} &:&
e_1 e_1 = e_3  & e_1 e_2=e_1 & e_2 e_1=e_3 &  e_2 e_2=e_4 \\
({\bf B}_2(\alpha))_{4,11} &:&
e_1 e_1 = e_3  & e_1 e_2=(1- \alpha) e_1 & e_2 e_1=\alpha e_1 &  e_2 e_2=e_4 \\
({\bf B}_2(\alpha\neq0))_{4,12}(\lambda\in {\bf k}) &:&
e_1 e_1 = e_3  & e_1 e_2=(1- \alpha) e_1+ \lambda e_3+e_4 & e_2 e_1=\alpha e_1 &  e_2 e_2=e_4 \\
({\bf B}_2(0))_{4,12}(\lambda\in {\bf k}) &:&
e_1 e_1 = e_3  & e_1 e_2=e_1 & e_2 e_1=\lambda  e_3+e_4 &  e_2 e_2=e_4 \\
({\bf B}_2(\alpha\neq0))_{5,13} &:&
e_1 e_1 = e_3  & e_1 e_2=(1- \alpha) e_1+e_4 & e_2 e_1=\alpha e_1&  e_2 e_2=e_5 \\
({\bf B}_2(0))_{5,13} &:&
e_1 e_1 = e_3  & e_1 e_2= e_1 & e_2 e_1=e_4 &  e_2 e_2=e_5
\end{array}$$

\subsection{Bilinear pair ${\bf B}_3$}

At first, we calculate $H^2({\bf B}_3,{\bf k})$  and $T_s({\bf B}_3)$.
It is easy to see that
$$\delta e_1^*=0, \delta e_2^*=\Delta_{12}- \Delta_{21},$$
and by Lemma \ref{l28},
$$Basis(H^2({\bf B}_3,{\bf k}))=
             \left\{ [\Delta_{11}], [\Delta_{21}], [\Delta_{22}] \right\}.$$

Since ${\rm rad}(g)=0$ then $T_s({\bf B}_3)=Grass_s(H^2({\bf B}_3,{\bf k}))$.

Now, we calculate the orbits of each space in $T_s({\bf B}_3)$ under the
action of ${\mathcal A}({\bf B}_3,g)$.

The action of ${\mathcal A}({\bf B}_3,g)$ on a subspace
$\langle a [\Delta_{11}]+ b[\Delta_{21}] +c [\Delta_{22}]  \rangle $
is the following:
$$\langle (a+bx+cx^2)[\Delta_{11}] + cxy [\Delta_{12}]+ (by+cxy)[\Delta_{21}]+ cy^2[\Delta_{22}] \rangle =
\langle (a+bx+cx^2)[\Delta_{11}] + (by+2cxy)[\Delta_{21}]+ cy^2[\Delta_{22}] \rangle.$$

Using the same techniques previously used, we get the following orbits:

\begin{enumerate}
    \item[$\star$]
     $T_1({\bf B}_3)=
    \{\langle [\Delta_{11}] \rangle \} \cup
      O( \langle [\Delta_{21}] \rangle ) \cup
      O( \langle [\Delta_{11}]+[\Delta_{22}] \rangle )\cup
      O( \langle [\Delta_{22}] \rangle ). $

    \item[$\star$]      $T_2({\bf B}_3)=
     O( \langle [\Delta_{21}], [\Delta_{22}] \rangle ) \cup
     O( \langle [\Delta_{11}]+ [\Delta_{22}], [\Delta_{21}]  \rangle)\cup
      O( \langle [\Delta_{11}], [\Delta_{21}]  \rangle)\cup  O( \langle [\Delta_{11}], [\Delta_{22}]  \rangle). $

    \item[$\star$]  $T_3({\bf B}_3)=H^2({\bf B}_3, {\bf k}).$
\end{enumerate}

At last, we compute each radical extension of $({\bf B}_3, g)$ using its
correspondence with an orbit of $T_s({\bf B}_3)$.
$$\begin{array}{ll llll}
({\bf B}_3)_{3,1} &:& e_1 e_1 = e_3  & e_1 e_2= e_2 & e_2 e_1=-e_2 & e_2 e_2=0\\
({\bf B}_3)_{3,2} &:& e_1 e_1 = 0  & e_1 e_2= e_2 & e_2 e_1=-e_2+e_3 & e_2 e_2=0\\
({\bf B}_3)_{3,3} &:& e_1 e_1 = e_3  & e_1 e_2= e_2 & e_2 e_1=-e_2 & e_2 e_2=e_3\\
({\bf B}_3)_{3,4} &:& e_1 e_1 = 0  & e_1 e_2= e_2 & e_2 e_1=-e_2 & e_2 e_2=e_3\\

({\bf B}_3)_{4,5} &:& e_1 e_1 = 0  & e_1 e_2= e_2 & e_2 e_1=-e_2+e_3 & e_2 e_2=e_4\\
({\bf B}_3)_{4,6} &:& e_1 e_1 = e_3  & e_1 e_2= e_2 & e_2 e_1=-e_2+e_4 & e_2 e_2=e_3\\
({\bf B}_3)_{4,7} &:& e_1 e_1 = e_3  & e_1 e_2= e_2 & e_2 e_1=-e_2+e_4 & e_2 e_2=0\\
({\bf B}_3)_{4,8} &:& e_1 e_1 = e_3  & e_1 e_2= e_2 & e_2 e_1=-e_2 & e_2 e_2=e_4\\

({\bf B}_3)_{5,9} &:& e_1 e_1 = e_3  & e_1 e_2= e_2 & e_2 e_1=-e_2+e_4 & e_2 e_2=e_5
\end{array}$$

\subsection{Bilinear pairs ${\bf C}(\alpha,\beta)$, $(\alpha, \beta) \in{\bf k}\times {\bf k_{\geq 0}}$}
Since the bilinear pairs ${\bf C}(0,\beta)$ and ${\bf C}(1,-\beta)^{op}$ are equivalent, we will  suppose now   $\alpha \neq 0.$

To begin with, we calculate $H^2({\bf C}(\alpha,\beta),{\bf k})$  and $T_s({\bf C}(\alpha, \beta))$.
It is easy to see that
$$\delta  e_1^*=(1- \alpha) \Delta_{12} + \alpha \Delta_{21}, \ \delta e_2^*=\Delta_{11}+\beta \Delta_{12}- \beta\Delta_{21}+\Delta_{22},$$
and by Lemma \ref{l28},
$$Basis(H^2({\bf C}(\alpha,\beta),{\bf k}))=
             \left\{ [\Delta_{12}], [\Delta_{22}] \right\}.$$

Since ${\rm rad}(g)=0$ then $T_s({\bf C}(\alpha, \beta))=Grass_s(H^2({\bf C}(\alpha,\beta),{\bf k}))$.

Moreover, we calculate the orbits of each space in $T_s({\bf C}(\alpha, \beta) )$ under the
action of ${\mathcal A}({\bf C}(\alpha, \beta),g)$.

The action of ${\mathcal A}({\bf C}(\alpha, \beta),g)$ on a vector space $\langle a[\Delta_{12}]+ b[\Delta_{22}] \rangle$
is the following:
$$\langle \gamma [\Delta_{12}] + b [\Delta_{22}]  \rangle, \mbox{
where }\gamma=\pm a\mbox{ if }\beta=0\mbox{ and }\gamma=a\mbox{ if }\beta \neq0.$$

As we did in the previous cases, we present the orbits:

\begin{enumerate}
 \item[$\star$]     $T_1({\bf C}(\alpha\neq0, 0))=   O(\langle [\Delta_{22}] \rangle ) \cup
     O(\{\langle [\Delta_{12}]+\lambda[\Delta_{22}] \rangle: \lambda \in {\bf k}_{\geq 0} \}).$
 \item[$\star$]     $T_1({\bf C}(\alpha\neq0, \beta))=   O( \langle [\Delta_{22}] \rangle) \cup
     O(\{\langle [\Delta_{12}]+\lambda[\Delta_{22}] \rangle: \lambda \in {\bf k} \})$ for $\beta\neq0$.

    \item[$\star$] $T_2({\bf C}(\alpha\neq0,\beta))=H^2({\bf C}(\alpha,\beta), {\bf k}).$

\end{enumerate}

To sum up, we compute each radical extension of $({\bf C}(\alpha,\beta), g)$ using its
correspondence with an orbit in $T_s({\bf C}(\alpha,\beta))$.
$$\begin{array}{ll llll}

({\bf C}(\alpha, \beta))_{3,1} &:&  e_1 e_1 = e_2 & e_1 e_2=(1-\alpha)e_1+ \beta e_2  & e_2 e_1=\alpha e_1 - \beta e_2 & e_2 e_2=e_2+e_3 \\

({\bf C}(\alpha\neq 0, 0))_{3,2}(\lambda\in{\bf k}_{\geq 0}) &:&  e_1 e_1 = e_2 & e_1 e_2=(1-\alpha)e_1 + e_3  & e_2 e_1=\alpha e_1 & e_2 e_2=e_2+\lambda e_3 \\

({\bf C}(0, 0))_{3,2}(\lambda\in{\bf k}_{\geq 0}) &:&  e_1 e_1 = e_2 & e_1 e_2= e_1 & e_2 e_1= e_3  & e_2 e_2=e_2+\lambda e_3 \\

({\bf C}(\alpha\neq 0, \beta\neq0))_{3,2}(\lambda\in{\bf k}) &:&  e_1 e_1 = e_2 & e_1 e_2=(1-\alpha)e_1+ \beta e_2+e_3 & e_2 e_1=\alpha e_1 - \beta e_2 & e_2 e_2=e_2+\lambda e_3 \\

({\bf C}(0, \beta\neq0))_{3,2}(\lambda\in{\bf k}) &:&  e_1 e_1 = e_2 & e_1 e_2=  e_1 + \beta e_2 & e_2 e_1= -\beta e_2+e_3 & e_2 e_2=e_2+\lambda e_3 \\

({\bf C}(\alpha\neq 0, \beta))_{4,3}(\lambda\in{\bf k}) &:&  e_1 e_1 = e_2 & e_1 e_2=(1-\alpha)e_1+ \beta e_2+e_3  & e_2 e_1=\alpha e_1 - \beta e_2  & e_2 e_2=e_2+e_4 \\

({\bf C}(0, \beta))_{4,3}(\lambda\in{\bf k}) &:&  e_1 e_1 = e_2 & e_1 e_2= e_1 + \beta e_2  & e_2 e_1=- \beta e_2 +e_3  & e_2 e_2=e_2+e_4
\end{array}$$

\subsection{Bilinear pairs ${\bf D}_1(\alpha,0)$, $(\alpha, 0)\in \mathcal{U}$}

Firstly, we calculate $H^2({\bf D}_1(\alpha,0),{\bf k})$  and $T_s({\bf D}_1(\alpha,0))$.
It is easy to see that
$$\delta  e_1^*=\Delta_{11}+(1- \alpha) \Delta_{12} + \alpha \Delta_{21}, \delta  e_2^*=0$$
and by Lemma \ref{l28},
$$Basis(H^2({\bf D}_1(\alpha,0),{\bf k}))=
             \left\{ [\Delta_{12}], [\Delta_{21}], [\Delta_{22}] \right\}.$$

Since ${\rm rad}(g)=0$ then $T_s({\bf D}_1(\alpha,0))=Grass_s(H^2({\bf D}_1(\alpha,0),{\bf k}))$.

The orbit of a space $\langle a[\Delta_{12}] +  b[\Delta_{21}]+ c[\Delta_{22}] \rangle$ under the action of an automorphism in ${\mathcal A}({\bf D}_1(\alpha,0),g)$ is the following set:
$$\{\langle a[\Delta_{12}]+  b[\Delta_{21}]+ c [\Delta_{22}]\rangle, \langle a[\Delta_{12}]+  b[\Delta_{21}]+ \gamma [\Delta_{22}] \rangle\}, \mbox{ where }  \gamma=(a+b-c) \mbox{ if }\alpha=\frac{1}{2} \mbox{ or }\gamma=c \mbox{ otherwise. }$$

The orbits are the following:
\begin{enumerate}
    \item[$\star$] $T_1({\bf D}_1(\alpha\neq \frac{1}{2},0))=
    \{ \langle [\Delta_{22}] \rangle \} \cup \{\langle [\Delta_{12}]+ \lambda [\Delta_{22}]  \rangle: \lambda \in {\bf k}\} \cup \{\langle [\Delta_{21}]+ \lambda [\Delta_{12}]+ \mu [\Delta_{22}]  \rangle: \lambda, \mu \in {\bf k}\} .$

   \item[$\star$] $T_1({\bf D}_1( \frac{1}{2},0))=
    \{ \langle [\Delta_{22}] \rangle \} \cup
   O( \{ \langle [\Delta_{12}]+\lambda[\Delta_{22}] \rangle: (\lambda, 0)\in \mathcal U\}) \cup
   O( \{ \langle [\Delta_{21}]+\lambda[\Delta_{12}]+\mu[\Delta_{22}] \rangle:(\mu, \lambda) \in \mathcal U  \}).$

    \item[$\star$]  $T_2({\bf D}_1(\alpha\neq \frac{1}{2},0))=$
    $ \{\langle [\Delta_{21}], [\Delta_{22}] \rangle \} \cup \{ \langle [\Delta_{22}], [\Delta_{12}]+\lambda[\Delta_{21}] \rangle:   \lambda \in {\bf k}\} \cup
    \{ \langle [\Delta_{21}]+\lambda[\Delta_{22}], [\Delta_{12}]+\mu[\Delta_{22}] \rangle:   \lambda,\mu \in {\bf k}\}. $
   \item[$\star$]  $T_2({\bf D}_1(\frac{1}{2},0))=$  $ O(\langle [\Delta_{21}], [\Delta_{22}] \rangle ) \cup
    O(\{ \langle [\Delta_{22}], [\Delta_{12}]+\lambda[\Delta_{21}] \rangle:   \lambda \in {\bf k}\}) \cup O(\{ \langle [\Delta_{21}]+\lambda[\Delta_{22}], [\Delta_{12}]+\mu[\Delta_{22}] \rangle:   (\lambda, 0), (\mu, 0) \in \mathcal U \}). $

    \item[$\star$]  $T_3({\bf D}_1(\alpha,0))=H^2({\bf D}_1(\alpha,0), {\bf k})$.
\end{enumerate}

Finally, we compute each radical extension  of $({\bf D}_1(\alpha,0), g)$ using its
correspondence with an orbit in $T_s({\bf D}_1(\alpha,0))$.
$$\begin{array}{ll llll}
({\bf D}_1(\alpha,0))_{3,1} &:&
e_1 e_1 = e_1 & e_1 e_2=(1-\alpha)e_1  & e_2 e_1=\alpha e_1 &e_2 e_2=e_3 \\
({\bf D}_1(\alpha\neq\frac{1}{2},0))_{3,2}(\lambda \in {\bf k}) &:&
e_1 e_1 = e_1 & e_1 e_2=(1-\alpha)e_1+e_3  & e_2 e_1=\alpha e_1 &e_2 e_2=\lambda e_3 \\
({\bf D}_1(\frac{1}{2},0)_{3,2}((\lambda, 0) \in \mathcal U) &:&
e_1 e_1 = e_1 & e_1 e_2=\frac{1}{2}e_1+e_3  & e_2 e_1=\frac{1}{2} e_1 &e_2 e_2=\lambda e_3 \\
({\bf D}_1(\alpha\neq\frac{1}{2},0))_{3,3}(\lambda, \mu \in {\bf k}) &:&
e_1 e_1 = e_1 & e_1 e_2=(1-\alpha)e_1+\lambda e_3  & e_2 e_1=\alpha e_1+e_3 &e_2 e_2=\mu e_3 \\
({\bf D}_1(\frac{1}{2},0))_{3,3}((\mu, \lambda) \in \mathcal U) &:&
e_1 e_1 = e_1 & e_1 e_2=\frac{1}{2}e_1+\lambda e_3  & e_2 e_1=\frac{1}{2} e_1+e_3 &e_2 e_2=\mu e_3 \\
({\bf D}_1(\alpha,0))_{4,5} &:&
e_1 e_1 = e_1 & e_1 e_2=(1-\alpha)e_1  & e_2 e_1=\alpha e_1+e_3 &e_2 e_2=e_4 \\
({\bf D}_1(\alpha,0))_{4,6}(\lambda \in {\bf k}) &:&
e_1 e_1 = e_1 & e_1 e_2=(1-\alpha)e_1 +e_4 & e_2 e_1=\alpha e_1+\lambda e_4 &e_2 e_2= e_3 \\
({\bf D}_1(\alpha\neq\frac{1}{2},0))_{4,7}(\lambda,\mu \in {\bf k}) &:&
e_1 e_1 = e_1 & e_1 e_2=(1-\alpha)e_1 +e_4 & e_2 e_1=\alpha e_1+ e_3 &e_2 e_2= \lambda e_3+\mu e_4 \\
({\bf D}_1(\frac{1}{2},0))_{4,7}((\lambda, 0), (\mu,0) \in \mathcal U) &:&
e_1 e_1 = e_1 & e_1 e_2=\frac{1}{2}e_1 +e_4 & e_2 e_1=\frac{1}{2} e_1+ e_3 &e_2 e_2= \lambda e_3+\mu e_4 \\
({\bf D}_1(\alpha,0))_{5,9} &:&
e_1 e_1 = e_1 & e_1 e_2=(1-\alpha)e_1+e_3  & e_2 e_1=\alpha e_1+e_4 &e_2 e_2=e_5
\end{array}$$

 \subsection{Bilinear pairs ${\bf D}_1(\alpha,\beta \neq 0)$, $(\alpha, \beta\neq 0)\in \mathcal{U}$}
To start with, we calculate $H^2({\bf D}_1(\alpha,\beta),{\bf k})$  and $T_s({\bf D}_1(\alpha,\beta))$.
It is easy to see that
$$\delta  e_1^*=\Delta_{11}+(1- \alpha) \Delta_{12} + \alpha \Delta_{21}, \ \delta e_2^*= \beta \Delta_{12}-\beta\Delta_{21},$$
and by Lemma \ref{l28},
$$Basis(H^2({\bf D}_1(\alpha,0),{\bf k}))=
             \left\{ [\Delta_{21}], [\Delta_{22}] \right\}.$$

Since ${\rm rad}(g)=0$ then $T_s({\bf D}_1(\alpha,\beta))=Grass_s(H^2({\bf D}_1(\alpha,\beta),{\bf k}))$.

Moreover, we calculate the orbits of each space in $T_s({\bf D}_1(\alpha,\beta))$ under the
action of ${\mathcal A}({\bf D}_1(\alpha,\beta),g)$.

Since ${\mathcal A}({\bf D}_1(\alpha,\beta),g)$ depends on the parameters, we consider two cases here. The orbit of a subspace
$\langle a [\Delta_{21}]+   b[\Delta_{22}]   \rangle$
by the action of ${\mathcal A}({\bf D}_1(\alpha,\beta),g)$ is
$\{\langle a[\Delta_{21}]+ b[\Delta_{22}]\rangle,\langle a[\Delta_{21}]+ (a-b)[\Delta_{22}]\rangle\}$ if $\beta=2\alpha-1$ or trivial otherwise. Thus, we obtain:

\begin{enumerate}
    \item[$\star$] $T_1({\bf D}_1(\alpha,2\alpha-1\neq0))= \{ \langle [\Delta_{22}] \rangle\} \cup O(\{ \langle [\Delta_{21}]+\lambda[\Delta_{22}]\rangle: (\lambda, 0)\in \mathcal U\}).$
    \item[$\star$] $T_1({\bf D}_1(\alpha,\beta\neq0))= \{ \langle [\Delta_{22}] \rangle\} \cup \{ \langle [\Delta_{21}]+\lambda[\Delta_{22}] \rangle: \lambda \in{\bf k}\} \textrm{ for } \beta\neq2\alpha-1 \textrm{ and } \beta\neq0$.

    \item[$\star$]  $T_2({\bf D}_1(\alpha,\beta\neq0))=H^2({\bf D}_1(\alpha,\beta), {\bf k})$.

\end{enumerate}

To sum up, we compute each radical extension of $({\bf D}_1(\alpha,\beta), g)$ using its
correspondence with an orbit in $T_s({\bf D}_1(\alpha,\beta))$.
{\small
$$\begin{array}{ll llll}
({\bf D}_1(\alpha,\beta\neq0))_{3,1}  &:&
e_1 e_1 = e_1 & e_1 e_2=(1-\alpha)e_1+ \beta e_2  & e_2 e_1=\alpha e_1 - \beta e_2 & e_2 e_2=e_3\\
({\bf D}_1(\alpha,2\alpha-1\neq0))_{3,4}( (\lambda, 0)\in \mathcal U)  &:&
e_1 e_1 = e_1 & e_1 e_2=(1-\alpha)e_1+ (2\alpha-1) e_2  & e_2 e_1=\alpha e_1 - (2\alpha-1) e_2+ e_3 & e_2 e_2=\lambda e_3\\
({\bf D}_1(\alpha,\beta\neq2\alpha-1, 0))_{3,4}(\lambda \in {\bf k}) &:&
e_1 e_1 = e_1 & e_1 e_2=(1-\alpha)e_1+ \beta e_2  & e_2 e_1=\alpha e_1 - \beta e_2+ e_3 & e_2 e_2=\lambda e_3\\
({\bf D}_1(\alpha,\beta\neq0))_{4,5}  &:&
e_1 e_1 = e_1 & e_1 e_2=(1-\alpha)e_1+ \beta e_2  & e_2 e_1=\alpha e_1 - \beta e_2+e_3 & e_2 e_2=e_4
\end{array}$$}

 \subsection{Bilinear pair ${\bf D}_2(0, 0)$}
At first, we calculate $H^2({\bf D}_2(0,0),{\bf k})$  and $T_s({\bf D}_2(0,0))$.
It is easy to see that $\delta  e_1^*=\Delta_{11},\delta  e_2^*=0.$ By Lemma \ref{l28},
$$Basis(H^2({\bf D}_2(0,0),{\bf k}))=
             \left\{ [\Delta_{12}], [\Delta_{21}], [\Delta_{22}] \right\}.$$

Also, it can be proved that $T_s({\bf D}_2(0,0))=Grass_s(H^2({\bf D}_2(0, 0),{\bf k}))$.

Now, we calculate the orbits of each space in $T_s({\bf D}_2(0, 0))$ under the
action of ${\mathcal A}({\bf D}_2(0,0),g)$.

The action of ${\mathcal A}({\bf D}_2(0,0),g)$ on a a subspace
$\langle  a[\Delta_{12}]+b[\Delta_{21}]+ c[\Delta_{22}] \rangle$ is
$$\langle ax[\Delta_{12}]+ bx [\Delta_{21}]+cx^2[\Delta_{22}]\rangle.$$

The orbits of the $T_s$ sets are presented below:
\begin{enumerate}
    \item[$\star$] $T_1({\bf D}_2(0,0))=$ $
    \{ \langle [\Delta_{22}] \rangle\} \cup
    O( \langle [\Delta_{21}]+[\Delta_{22}] \rangle) \cup
    \{ \langle [\Delta_{21}] \rangle\} \cup
    O(\{ \langle [\Delta_{12}]+ \lambda [\Delta_{21}]+ [\Delta_{22}]\rangle: \lambda \in {\bf k}\}) \cup$

    $\cup O(\{ \langle [\Delta_{12}] + \lambda [\Delta_{21}] \rangle: \lambda \in {\bf k}\}).$

    \item[$\star$]  $T_2({\bf D}_2(0, 0))=$ $
    O( \langle [\Delta_{21}], [\Delta_{22}] \rangle) \cup
    O(\langle [\Delta_{12}], [\Delta_{21}]+[\Delta_{22}] \rangle) \cup
    O( \langle  [\Delta_{12}], [\Delta_{21}] \rangle) \cup$

    $O(\{ \langle [\Delta_{22}], [\Delta_{12}]+ \lambda [\Delta_{21}] \rangle: \lambda \in {\bf k}\})\cup O(\{ \langle [\Delta_{21}]+ \lambda[\Delta_{22}], [\Delta_{12}]+[\Delta_{22}] \rangle: \lambda \in {\bf k}\}).$

    \item[$\star$]  $T_3({\bf D}_2(0, 0))=H^2({\bf D}_2(0,0), {\bf k}).$

\end{enumerate}

At last, we compute each radical extension of $({\bf D}_2(0,0), g)$ using its
correspondence with an orbit in $T_s({\bf D}_2(0,0))$.
$$\begin{array}{ll llll}
({\bf D}_2(0,0))_{3,1} &:& e_1 e_1 = e_1 & e_1 e_2=0 & e_2 e_1= 0 & e_2 e_2=e_3\\
({\bf D}_2(0,0))_{3,2} &:& e_1 e_1 = e_1 & e_1 e_2=0 & e_2 e_1= e_3 & e_2 e_2=e_3\\
({\bf D}_2(0,0))_{3,3} &:& e_1 e_1 = e_1 & e_1 e_2=0 & e_2 e_1= e_3 & e_2 e_2=0\\
({\bf D}_2(0,0))_{3,4}(\lambda\in {\bf k}) &:& e_1 e_1 = e_1 & e_1 e_2=e_3 & e_2 e_1=\lambda e_3 & e_2 e_2=e_3\\
({\bf D}_2(0,0))_{3,5}(\lambda\in {\bf k}) &:& e_1 e_1 = e_1 & e_1 e_2=e_3 & e_2 e_1=\lambda e_3 & e_2 e_2=0\\
({\bf D}_2(0,0))_{4,8} &:& e_1 e_1 = e_1 & e_1 e_2=0 & e_2 e_1= e_3 & e_2 e_2=e_4 \\
({\bf D}_2(0,0))_{4,9} &:& e_1 e_1 = e_1 & e_1 e_2=e_3 & e_2 e_1= e_4 & e_2 e_2=e_4 \\
({\bf D}_2(0,0))_{4,10} &:& e_1 e_1 = e_1 & e_1 e_2=e_3 & e_2 e_1= e_4 & e_2 e_2=0 \\
({\bf D}_2(0,0))_{4,11}(\lambda\in {\bf k}) &:& e_1 e_1 = e_1 & e_1 e_2=e_4 & e_2 e_1= \lambda e_4 & e_2 e_2=e_3 \\
({\bf D}_2(0,0))_{4,12}(\lambda\in {\bf k}) &:& e_1 e_1 = e_1 & e_1 e_2=e_4 & e_2 e_1= e_3 & e_2 e_2=\lambda e_3+e_4 \\
({\bf D}_2(0,0))_{5,14} &:& e_1 e_1 = e_1 & e_1 e_2=e_3 & e_2 e_1= e_4 & e_2 e_2=e_5
\end{array}$$

 \subsection{Bilinear pairs ${\bf D}_2(\alpha, \beta)$, $(0,0)\neq(\alpha,\beta)\in{\bf k}^2\setminus \mathcal{T}$.}
  
  Since the bilinear pairs ${\bf D}_2(\alpha, 0)$ and ${\bf D}_2(0, \alpha)^{op}$ are equivalent, we will  suppose now $\beta \neq 0$.

It is easy to see that
$$\delta  e_1^*=\Delta_{11}, \ \delta e_2^*=\alpha \Delta_{12}+ \beta \Delta_{21},$$
and by Lemma \ref{l28},
$$Basis(H^2({\bf D}_2(\alpha,\beta),{\bf k}))=
             \left\{ [\Delta_{12}],  [\Delta_{22}] \right\}.$$

Since ${\rm rad}(g)=0$ then $T_s({\bf D}_2(\alpha,\beta))=Grass_s(H^2({\bf D}_2(\alpha, \beta),{\bf k}))$.

Additionally, we calculate the orbits of each space in $T_s({\bf D}_2(\alpha, \beta))$ under the
action of ${\mathcal A}({\bf D}_2(\alpha,\beta),g)$.

The action of ${\mathcal A}({\bf D}_2(\alpha,\beta),g$ on a subspace
$\langle a[\Delta_{12}]+ c[\Delta_{22}]  \rangle$
is
$\langle ax[\Delta_{12}]+ cx^2[\Delta_{22}] \rangle.$

It is easy to see that
\begin{enumerate}
    \item[$\star$] $T_1({\bf D}_2(\alpha,\beta\neq0))=
    \{\langle [\Delta_{22}] \rangle\} \cup O(\langle  [\Delta_{12}] +[\Delta_{22}] \rangle) \cup \{\langle [\Delta_{12}] \rangle\}.$

    \item[$\star$] $T_2({\bf D}_2(\alpha, \beta\neq0))=H^2({\bf D}_2(\alpha,\beta), {\bf k}).$

\end{enumerate}

Lastly, we compute each radical extension of $({\bf D}_2(\alpha,\beta), g)$ using its
correspondence with an orbit in $T_s({\bf D}_2(\alpha,\beta))$.
$$\begin{array}{ll llll}
({\bf D}_2(\alpha,\beta))_{3,1} &:& e_1 e_1 = e_1 & e_1 e_2=\alpha e_2  & e_2 e_1= \beta e_2 & e_2 e_2=e_3 \\
({\bf D}_2(\alpha,\beta\neq0))_{3,6} &:& e_1 e_1 = e_1 & e_1 e_2=\alpha e_2 +e_3 & e_2 e_1= \beta e_2 & e_2 e_2=e_3 \\
({\bf D}_2(\alpha,0))_{3,2}  &:& e_1 e_1 = e_1 & e_1 e_2=\alpha e_2 & e_2 e_1= e_3 & e_2 e_2=e_3 \\
({\bf D}_2(\alpha,\beta\neq0))_{3,7} &:& e_1 e_1 = e_1 & e_1 e_2=\alpha e_2+e_3 & e_2 e_1= \beta e_2  & e_2 e_2=0 \\
({\bf D}_2(\alpha,0))_{3,3} &:& e_1 e_1 = e_1 & e_1 e_2=\alpha e_2 & e_2 e_1= e_3  & e_2 e_2=0 \\
({\bf D}_2(\alpha,\beta\neq0))_{4,13} &:& e_1 e_1 = e_1 & e_1 e_2=\alpha e_2+e_3 & e_2 e_1= \beta e_2  & e_2 e_2=e_4\\
({\bf D}_2(\alpha,0))_{4,8} &:& e_1 e_1 = e_1 & e_1 e_2=\alpha e_2 & e_2 e_1= e_3  & e_2 e_2=e_4
\end{array}$$

\subsection{Bilinear pair ${\bf D}_3(0,0)$}
Firstly, we calculate $H^2({\bf D}_3(0,0),{\bf k})$  and $T_s({\bf D}_3(0,0))$.
It is easy to see that
$$\delta  e_1^*=\Delta_{11}+ \Delta_{12} - \Delta_{21}, \delta  e_2^*=0$$
and by Lemma \ref{l28},
$$Basis(H^2({\bf D}_3(\alpha,\beta),{\bf k}))=
             \left\{ [\Delta_{12}],  [\Delta_{21}], [\Delta_{22}] \right\}.$$

Since ${\rm rad}(g)=0$ then $T_s({\bf D}_3(0,0))=Grass_s(H^2({\bf D}_3(0,0),{\bf k}))$.

Secondly,
we note that  the action of ${\mathcal A}({\bf D}_3(\alpha,\beta),g)$ is the identity.

\begin{enumerate}
    \item[$\star$] $T_1({\bf D}_3(0,0))=
    \{\langle [\Delta_{22}] \rangle\} \cup
    \{ \langle [\Delta_{12}]+\lambda [\Delta_{22}]\rangle: \lambda \in {\bf k} \} \cup
    \{ \langle [\Delta_{21}]+\lambda [\Delta_{12}]+\mu [\Delta_{22}]\rangle: \lambda, \mu \in {\bf k} \}.$

    \item[$\star$] $T_2({\bf D}_3(0,0))=$
    $\{\langle [\Delta_{21}], [\Delta_{22}] \rangle\} \cup \{ \langle [\Delta_{22}], [\Delta_{12}] + \lambda [\Delta_{21}]\rangle:\lambda \in {\bf k} \} \cup
    \{ \langle [\Delta_{21}]+\lambda [\Delta_{22}], [\Delta_{12}]+\mu [\Delta_{22}]\rangle:\lambda, \mu \in {\bf k} \}.$

    \item[$\star$]  $T_3({\bf D}_3(0,0))=H^2({\bf D}_3(0,0), {\bf k})$.

\end{enumerate}

Finally, we compute each radical extension of $({\bf D}_3(0,0), g)$ using its
correspondence with an orbit in $T_s({\bf D}_3(0,0))$.
$$\begin{array}{ll llll}
({\bf D}_3(0,0))_{3,1} &:&
e_1 e_1 = e_1 & e_1 e_2=e_1    & e_2 e_1=- e_1 &  e_2 e_2=e_3 \\
({\bf D}_3(0,0))_{3,2}(\lambda \in {\bf k}) &:&
e_1 e_1 = e_1 & e_1 e_2=e_1+e_3    & e_2 e_1=- e_1 &  e_2 e_2=\lambda e_3\\
({\bf D}_3(0,0))_{3,3}(\lambda, \mu\in {\bf k}) &:&
e_1 e_1 = e_1 & e_1 e_2=e_1+\lambda e_3    & e_2 e_1=- e_1+e_3 &  e_2 e_2=\mu e_3 \\
({\bf D}_3(0,0))_{4,5} &:&
e_1 e_1 = e_1 & e_1 e_2=e_1   & e_2 e_1=- e_1+e_3 &  e_2 e_2=e_4 \\
({\bf D}_3(0,0))_{4,6}(\lambda \in {\bf k}) &:&
e_1 e_1 = e_1 & e_1 e_2=e_1+e_4    & e_2 e_1=- e_1+\lambda e_4 &  e_2 e_2=e_3 \\
({\bf D}_3(0,0))_{4,7}(\lambda, \mu\in {\bf k}) &:&
e_1 e_1 = e_1 & e_1 e_2=e_1 + e_4    & e_2 e_1=- e_1+e_3 &  e_2 e_2=\lambda e_3 +\mu e_4 \\
({\bf D}_3(0,0))_{5,9} &:&
e_1 e_1 = e_1 & e_1 e_2=e_1+e_3    & e_2 e_1=- e_1+e_4 &  e_2 e_2=e_5

\end{array}$$

 \subsection{Bilinear pairs ${\bf D}_3(0,\beta)$, $(0,\beta)\in{\bf k}^2\setminus (\mathcal{T}\cup \{(0,0)\})$}

To begin with, we calculate $H^2({\bf D}_3(0,\beta),{\bf k})$  and $T_s({\bf D}_3(0,\beta))$ for $\beta\neq0$.
It is easy to see that
$$\delta  e_1^*=\Delta_{11}+ \Delta_{12} - \Delta_{21}, \ \delta e_2^*= \beta\Delta_{21},$$
and by Lemma \ref{l28},
$$Basis(H^2({\bf D}_3(0,\beta),{\bf k}))=
             \left\{  [\Delta_{12}], [\Delta_{22}] \right\}.$$

Since ${\rm rad}(g)=0$ then $T_s({\bf D}_3(0,\beta))=Grass_s(H^2({\bf D}_3(0,\beta),{\bf k}))$.

Since the action of ${\mathcal A}({\bf D}_3(0,\beta),g)$ is the identity, we get the following orbits

\begin{enumerate}
    \item[$\star$] $T_1({\bf D}_3(0,\beta\neq0))=
     \{\langle [\Delta_{22}] \rangle\} \cup \{ \langle [\Delta_{12}]+\lambda[\Delta_{22}] \rangle: \lambda \in {\bf k}\}$.

    \item[$\star$]  $T_2({\bf D}_3(0,\beta\neq0))=H^2({\bf D}_3(0,\beta), {\bf k})$.

\end{enumerate}

To sum up, we compute each radical extension of $({\bf D}_3(0,\beta), g)$ using its
correspondence with an orbit in $T_s({\bf D}_3(0,\beta))$.
$$\begin{array}{ll llll}
({\bf D}_3(0,\beta\neq0))_{3,1} &:&
e_1 e_1 = e_1 & e_1 e_2= e_1    & e_2 e_1= -e_1+\beta e_2   &  e_2 e_2=e_3 \\
({\bf D}_3(0,\beta\neq0))_{3,2}(\lambda\in {\bf k}) &:&
e_1 e_1 = e_1 & e_1 e_2= e_1 +e_3   & e_2 e_1= -e_1+\beta e_2   &  e_2 e_2=\lambda e_3 \\
({\bf D}_3(0,\beta\neq0))_{4,8} &:&
e_1 e_1 = e_1 & e_1 e_2= e_1 +e_3     & e_2 e_1= -e_1+\beta e_2       &  e_2 e_2=e_4
\end{array}$$


 \subsection{Bilinear pairs ${\bf D}_3(\alpha \neq0,\beta)$, $(\alpha,\beta)\in{\bf k}^2\setminus \mathcal{T}$}

To begin with, we calculate $H^2({\bf D}_3(\alpha,\beta),{\bf k})$  and $T_s({\bf D}_3(\alpha,\beta))$.
It is easy to see that
$$\delta  e_1^*=\Delta_{11}+ \Delta_{12} - \Delta_{21}, \ \delta e_2^*= \alpha \Delta_{12}+\beta\Delta_{21},$$
and by Lemma \ref{l28},
$$Basis(H^2({\bf D}_3(\alpha,\beta),{\bf k}))=
             \left\{  [\Delta_{21}], [\Delta_{22}] \right\}.$$

Since ${\rm rad}(g)=0$ then $T_s({\bf D}_3(\alpha,\beta))=Grass_s(H^2({\bf D}_3(\alpha,\beta),{\bf k}))$.

Also, the action of ${\mathcal A}({\bf D}_3(\alpha,\beta),g)$ is the identity action.

\begin{enumerate}
    \item[$\star$] $T_1({\bf D}_3(\alpha\neq0,\beta))=
     \{\langle [\Delta_{22}] \rangle\} \cup \{ \langle [\Delta_{21}]+\lambda[\Delta_{22}] \rangle: \lambda \in {\bf k}\}$.

    \item[$\star$]  $T_2({\bf D}_3(\alpha\neq0,\beta))=H^2({\bf D}_3(\alpha,\beta), {\bf k})$.

\end{enumerate}

To sum up, we compute each radical extension of $({\bf D}_3(\alpha,\beta), g)$ using its
correspondence with an orbit in $T_s({\bf D}_3(\alpha,\beta))$.
$$\begin{array}{ll llll}
({\bf D}_3(\alpha\neq0,\beta))_{3,1} &:&
e_1 e_1 = e_1 & e_1 e_2=e_1+ \alpha e_2   & e_2 e_1=- e_1 + \beta e_2 &  e_2 e_2=e_3 \\
({\bf D}_3(\alpha\neq0,\beta))_{3,4}(\lambda\in {\bf k}) &:&
e_1 e_1 = e_1 & e_1 e_2=e_1+ \alpha e_2   & e_2 e_1=- e_1 + \beta e_2+ e_3 &  e_2 e_2=\lambda e_3 \\
({\bf D}_3(\alpha\neq0,\beta))_{4,5} &:&
e_1 e_1 = e_1 & e_1 e_2=e_1+ \alpha e_2   & e_2 e_1=- e_1 + \beta e_2+e_3 &  e_2 e_2=e_4
\end{array}$$


\subsection{Bilinear pairs ${\bf E}_1(\alpha,\beta,\gamma,\delta)$, $(\alpha,\beta,\gamma,\delta)\in \mathcal{V}$}
At first, we calculate $H^2({\bf E}_1(\alpha,\beta,\gamma,\delta), {\bf k})$ and $T_s({\bf E}_1(\alpha,\beta,\gamma,\delta)).$

It is easy to see  that
$$\delta e_1^*=\Delta_{11}+\alpha\Delta_{12}+\gamma  \Delta_{21}, \ \delta e_2^*=\beta \Delta_{12}+\delta \Delta_{21}+\Delta_{22},$$

and by Lemma \ref{l28},
$$Basis(H^2({\bf E}_1(\alpha,\beta,\gamma,\delta),{\bf k}))=
             \{ [\Delta_{12}], [\Delta_{21}]  \}.$$

Since ${\rm rad}(g)=0$ then $T_s({\bf E}_1(\alpha,\beta,\gamma,\delta))=Grass_s(H^2({\bf E}_1(\alpha,\beta,\gamma,\delta),{\bf k}))$.

The orbit of a subspace
$\langle a[\Delta_{12}]+ b[\Delta_{21}]  \rangle $
by the action of ${\mathcal A}({\bf E}_1(\alpha,\beta,\gamma,\delta),g)$ is the following set:
$$\{ \langle a[\Delta_{12}]+ b[\Delta_{21}] \rangle,  \langle b[\Delta_{12}]+ a[\Delta_{21}] \rangle\},$$
if $\alpha=\delta$ and $\beta=\gamma,$ otherwise it is the identity action.

It is easy to see that
\begin{enumerate}
 \item[$\star$] $T_1({\bf E}_1(\alpha,\beta,\gamma,\delta))=
    \{\langle [\Delta_{21}] \rangle\} \cup \{ \langle \lambda[\Delta_{21}]+[\Delta_{12}] \rangle:\lambda \in {\bf k}\}$ for $\gamma\neq\beta \textrm{ or } \delta\neq\alpha$.

 \item[$\star$] $T_1({\bf E}_1(\alpha,\beta,\beta,\alpha))=
    O( \langle [\Delta_{21}] \rangle) \cup O(\{ \langle \lambda[\Delta_{21}]+[\Delta_{12}] \rangle: \lambda \in {\bf k}^*_{>1}\cup \{1\}\})$.

    \item[$\star$]  $T_2({\bf E}_1(\alpha,\beta,\gamma,\delta))=H^2({\bf E}_1(\alpha,\beta,\gamma,\delta), {\bf k})$.

\end{enumerate}

At last, we compute each radical extension of $({\bf E}_1(\alpha,\beta,\gamma), g)$ using its
correspondence with an orbit in $T_s({\bf E}_1(\alpha,\beta,\gamma)):$
$$\begin{array}{ll llll}
({\bf E}_1(\alpha, \beta, \gamma,\delta))_{3,1} &:&  e_1 e_1 = e_1 & e_1 e_2=\alpha e_1+ \beta e_2 & e_2 e_1= \gamma  e_1 + \delta  e_2 +e_3& e_2 e_2=e_2 \\
({\bf E}_1(\alpha, \beta, \gamma,\delta))_{3,2}(\lambda \in {\bf k})\\
(\gamma\neq\beta \textrm{ or } \delta\neq\alpha) &:&  e_1 e_1 = e_1 & e_1 e_2=\alpha e_1+ \beta e_2+e_3 & e_2 e_1= \gamma  e_1 + \delta  e_2 +\lambda e_3& e_2 e_2=e_2 \\
({\bf E}_1(\alpha, \beta, \beta,\alpha))_{3,2}(\lambda \in {\bf k}^*_{>1}\cup \{1\}) &:&  e_1 e_1 = e_1 & e_1 e_2=\alpha e_1+ \beta e_2+e_3 & e_2 e_1= \beta  e_1 + \alpha  e_2 +\lambda e_3& e_2 e_2=e_2 \\
({\bf E}_1(\alpha, \beta, \gamma,\delta))_{4,3} &:&  e_1 e_1 = e_1 & e_1 e_2=\alpha e_1+ \beta e_2+e_3 & e_2 e_1= \gamma  e_1 + \delta  e_2 + e_4 & e_2 e_2=e_2
\end{array}$$

\subsection{Bilinear pairs ${\bf E}_2(\alpha,\beta,\gamma)$,  $(\alpha,\beta,\gamma)\in {\bf k}^3\setminus {\bf k}\times\mathcal{T}$}

To begin with, we calculate $H^2({\bf E}_2(\alpha,\beta,\gamma),{\bf k})$ and $T_s({\bf E}_2(\alpha,\beta,\gamma))$.
It is easy to see  that
$$\delta e_1^*=\Delta_{11}+(1-\alpha)\Delta_{12}+\alpha \Delta_{21}, \ \delta e_2^*=\beta \Delta_{12}+\gamma \Delta_{21}+\Delta_{22},$$

and by Lemma \ref{l28},
$$Basis(H^2({\bf E}_2(\alpha,\beta,\gamma),{\bf k}))=
             \left\{ [\Delta_{12}], [\Delta_{21}]  \right\}.$$

Since ${\rm rad}(g)=0$ then $T_s({\bf E}_2(\alpha,\beta,\gamma))=Grass_s(H^2({\bf E}_2(\alpha,\beta,\gamma),{\bf k}))$.

Additionally, we note that the action of ${\mathcal A}({\bf E}_2(\alpha,\beta,\gamma),g)$ is the identity and therefore:
\begin{enumerate}
     \item[$\star$] $T_1({\bf E}_2(\alpha,\beta,\gamma))=
    \{\langle [\Delta_{21}] \rangle\} \cup \{ \langle \lambda[\Delta_{21}]+[\Delta_{12}] \rangle:\lambda\in{\bf k}\}$.

    \item[$\star$]  $T_2({\bf E}_2(\alpha,\beta,\gamma))=H^2({\bf E}_2(\alpha,\beta,\gamma), {\bf k})$.

\end{enumerate}

Finally, we compute each radical extension of $({\bf E}_2(\alpha,\beta,\gamma), g)$ using its
correspondence with an orbit in $T_s({\bf E}_2(\alpha,\beta,\gamma)):$
$$\begin{array}{ll llll}
({\bf E}_2(\alpha, \beta, \gamma))_{3,1} &:&  e_1 e_1 = e_1 & e_1 e_2=(1-\alpha)e_1+ \beta e_2 & e_2 e_1= \alpha e_1 + \gamma e_2 +e_3& e_2 e_2=e_2 \\
({\bf E}_2(\alpha, \beta, \gamma))_{3,2}(\lambda\in {\bf k}) &:&  e_1 e_1 = e_1 & e_1 e_2=(1-\alpha)e_1+ \beta e_2+e_3 & e_2 e_1= \alpha e_1 + \gamma e_2 +\lambda e_3& e_2 e_2=e_2 \\
({\bf E}_2(\alpha, \beta, \gamma))_{4,3} &:&  e_1 e_1 = e_1 & e_1 e_2=(1-\alpha)e_1+ \beta e_2+e_3 & e_2 e_1= \alpha e_1 + \gamma e_2+e_4 & e_2 e_2=e_2

\end{array}$$

 \subsection{Bilinear pairs ${\bf E}_3(\alpha,\beta,\gamma)$, $(\alpha,\beta,\gamma)\in {\bf k}^2\times{\bf k}^*_{>1}$}

We calculate $H^2({\bf E}_3(\alpha,\beta,\gamma),{\bf k})$ and $T_s({\bf E}_3(\alpha,\beta,\gamma))$.
It is easy to see  that
$$\delta e_1^*=\Delta_{11}+(1-\alpha)\gamma \Delta_{12}+\alpha\gamma \Delta_{21}, \ \delta e_2^*=\frac{\beta}{\gamma} \Delta_{12}+\frac{1-\beta}{\gamma} \Delta_{21}+\Delta_{22}.$$

and by Lemma \ref{l28},
$$Basis(H^2({\bf E}_3(\alpha,\beta,\gamma),{\bf k}))=
             \left\{ [\Delta_{12}], [\Delta_{21}]  \right\}.$$

Since ${\rm rad}(g)=0$ then $T_s({\bf E}_3(\alpha,\beta,\gamma))=Grass_s(H^2({\bf E}_3(\alpha,\beta,\gamma),{\bf k}))$.

The orbit of a subspace
$\langle a[\Delta_{12}]+ b[\Delta_{21}]  \rangle $
by the action of ${\mathcal A}({\bf E}_3(\alpha,\beta,\gamma),g)$ is the following set:
$$\{ \langle a[\Delta_{12}]+ b[\Delta_{21}] \rangle,  \langle b[\Delta_{12}]+ a[\Delta_{21}] \rangle\},$$
if $\gamma=-1$ and $\alpha=\beta,$
otherwise it is the identity action. Thus, we have distinguished two cases for $T_1({\bf E}_3(\alpha,\beta,\gamma))$

\begin{enumerate}
     \item[$\star$] $T_1({\bf E}_3(\alpha, \beta, \gamma)=
    \{\langle [\Delta_{21}] \rangle\} \cup \{ \langle \lambda[\Delta_{21}]+[\Delta_{12}] \rangle:\lambda\in{\bf k}\}$ for $\beta\neq \alpha \textrm{ or } \gamma\neq -1$.
     \item[$\star$] $T_1({\bf E}_3(\alpha,\alpha,-1))=
    O(\langle [\Delta_{21}] \rangle) \cup O( \{\langle \lambda[\Delta_{21}]+[\Delta_{12}] \rangle: \lambda \in {\bf k}^*_{>1}\cup \{1\}\})$.

    \item[$\star$]  $T_2({\bf E}_3(\alpha,\beta,\gamma))=H^2({\bf E}_3(\alpha,\beta,\gamma), {\bf k})$.

\end{enumerate}

Now, we compute each radical extension of $({\bf E}_3(\alpha,\beta,\gamma), g)$ using its
correspondence with an orbit in $T_s({\bf E}_3(\alpha,\beta,\gamma)):$
$$\begin{array}{ll llll}
({\bf E}_3(\alpha, \beta, \gamma))_{3,1} &:&  e_1 e_1 = e_1 & e_1 e_2=(1-\alpha)\gamma e_1+ \frac{\beta}{\gamma} e_2  & e_2 e_1=\alpha \gamma e_1 + \frac{1-\beta}{\gamma}e_2 +e_3&
e_2 e_2=e_2 \\
({\bf E}_3(\alpha, \beta, \gamma))_{3,2}(\lambda\in {\bf k}) \\
(\beta\neq \alpha \textrm{ or } \gamma\neq -1) &:&  e_1 e_1 = e_1 & e_1 e_2=(1-\alpha)\gamma e_1+ \frac{\beta}{\gamma} e_2  +e_3& e_2 e_1=\alpha \gamma e_1 + \frac{1-\beta}{\gamma}e_2+\lambda e_3 &
e_2 e_2=e_2\\
({\bf E}_3(\alpha, \alpha, -1))_{3,2} &&&\\
(\lambda\in {\bf k}^*_{>1}\cup \{1\}) &:&  e_1 e_1 = e_1 & e_1 e_2=-(1-\alpha) e_1-\alpha e_2  +e_3& e_2 e_1=-\alpha e_1  -(1-\alpha)e_2+\lambda e_3 &
e_2 e_2=e_2 \\
({\bf E}_3(\alpha, \beta, \gamma))_{4,3} &:&  e_1 e_1 = e_1 & e_1 e_2=(1-\alpha)\gamma e_1+ \frac{\beta}{\gamma} e_2+e_3  & e_2 e_1=\alpha \gamma e_1 + \frac{1-\beta}{\gamma}e_2 +e_4&
e_2 e_2=e_2
\end{array}$$

 \subsection{Bilinear pair ${\bf E}_4$}

At first, we calculate $H^2({\bf E}_4,{\bf k})$ and $T_s({\bf E}_4)$.
It is easy to see  that
$$\delta e_1^*=\Delta_{11}+\Delta_{12}, \ \delta e_2^*=\Delta_{12}+\Delta_{22}.$$

and by Lemma \ref{l28},
$$Basis(H^2({\bf E}_4,{\bf k}))=
             \left\{ [\Delta_{12}], [\Delta_{21}]  \right\}.$$

Since ${\rm rad}(g)=0$ then $T_s({\bf E}_4)=Grass_s(H^2({\bf E}_4,{\bf k}))$.

Moreover, we note that the action of ${\mathcal A}({\bf E}_4,g)$ is the identity. Therefore we have the following orbits:

\begin{enumerate}
     \item[$\star$] $T_1({\bf E}_4)=
    \{\langle [\Delta_{21}] \rangle\} \cup\{ \langle \lambda[\Delta_{21}]+[\Delta_{12}] \rangle: \lambda \in {\bf k}\}$.

    \item[$\star$]  $T_2({\bf E}_4)=H^2({\bf E}_4, {\bf k})$.

\end{enumerate}

At last, we compute each radical extension of $({\bf E}_4, g)$ using its
correspondence with an orbit in $T_s({\bf E}_4):$
$$\begin{array}{ll llll}
({\bf E}_4)_{3,1} &:&  e_1 e_1 = e_1 & e_1 e_2=e_1+ e_2  &  e_2 e_1=e_3  & e_2 e_2=e_2 \\
({\bf E}_4)_{3,2}(\lambda\in {\bf k}) &:&  e_1 e_1 = e_1 & e_1 e_2=e_1+ e_2+e_3  &  e_2 e_1=\lambda e_3  & e_2 e_2=e_2 \\
({\bf E}_4)_{4,3} &:&  e_1 e_1 = e_1 & e_1 e_2=e_1+ e_2+e_3  &  e_2 e_1= e_4  & e_2 e_2=e_2
\end{array}$$

\subsection{Bilinear pairs ${\bf E}_5(\alpha)$, $\alpha\in{\bf k}$}

Firstly, we calculate $H^2({\bf E}_5(\alpha),{\bf k})$ and $T_s({\bf E}_5(\alpha))$.
It is easy to see that
$$\delta e_1^* = \Delta_{11}+(1-\alpha) \Delta_{12}+\alpha \Delta_{21}, \ \delta e_2^*= \alpha \Delta_{12}+(1-\alpha)\Delta_{21}+\Delta_{22},$$
and by Lemma \ref{l28},
$$Basis(H^2({\bf E}_5(\alpha),{\bf k}))=
             \left\{ [\Delta_{12}], [\Delta_{21}] \right\}.$$

Since ${\rm rad}(g)=0$ then $T_s({\bf E}_5(\alpha))=Grass_s(H^2({\bf E}_5(\alpha),{\bf k}))$.

Secondly, we calculate the orbits of each space in $T_s({\bf E}_5(\alpha))$ under the
action of ${\mathcal A}({\bf E}_5(\alpha),g)$.

The action of ${\mathcal A}({\bf E}_5(\alpha),g)$ on a vector space $\langle a[\Delta_{12}] +  b[\Delta_{21}]  \rangle$ is the following:

$\langle (a(1-y)x+b(1-x)y-(1-\alpha)(a+b)(1-x)x-\alpha(a+b)(1-y)y)[\Delta_{12}] + $
$$+(a(1-x)y+b(1-y)x-(1-\alpha)(a+b)(1-y)y-\alpha(a+b)(1-x)x)[\Delta_{21}] \rangle.$$

It can be proved that:
\begin{enumerate}
    \item[$\star$] $T_1({\bf E}_5(\alpha\neq 1/2))=O( \langle [\Delta_{21}]\rangle ) \cup O( \langle [\Delta_{12}]-[\Delta_{21}]\rangle ).$
        \item[$\star$] $T_1({\bf E}_5(1/2))=O( \langle [\Delta_{21}]\rangle ) \cup O(\langle [\Delta_{12}]-[\Delta_{21}]\rangle)\cup O( \langle [\Delta_{12}]+[\Delta_{21}]\rangle).$
    \item[$\star$] $T_2({\bf E}_5(\alpha))=H^2({\bf E}_5(\alpha), {\bf k})$.
\end{enumerate}

Finally, we compute each radical extension of $({\bf E}_5(\alpha), g)$ using its
correspondence with an orbit in $T_s({\bf E}_5(\alpha)):$
$$\begin{array}{ll llll}
({\bf E}_5(\alpha))_{3,1} &:& e_1 e_1 = e_1 & e_1 e_2=(1-\alpha)e_1+ \alpha e_2 & e_2 e_1=\alpha e_1+ (1-\alpha)e_2+e_3   & e_2 e_2=e_2 \\
({\bf E}_5(\alpha))_{3,2} &:& e_1 e_1 = e_1 & e_1 e_2=(1-\alpha)e_1+ \alpha e_2+e_3 & e_2 e_1=\alpha e_1+ (1-\alpha) e_2-e_3  & e_2 e_2=e_2 \\
({\bf E}_5(1/2))_{3,3} &:& e_1 e_1 = e_1 & e_1 e_2=\frac{1}{2}e_1+ \frac{1}{2} e_2+e_3 & e_2 e_1=\frac{1}{2} e_1+ \frac{1}{2} e_2+e_3  & e_2 e_2=e_2 \\
({\bf E}_5(\alpha))_{4,4} &:& e_1 e_1 = e_1 & e_1 e_2=(1-\alpha)e_1+ \alpha e_2+e_3 & e_2 e_1=\alpha e_1+ (1-\alpha) e_2+e_4  & e_2 e_2=e_2
\end{array}$$

\subsection{Bilinear pair $\mathfrak{N}_2$}

As before, we start calculating  $H^2(\mathfrak{N}_2,{\bf k})$  and $T_s(\mathfrak{N}_2)$. Since ${B_g(\mathfrak{N}_2 \times \mathfrak{N}_2, \mathbb {\bf k})}=\left\{0\right\}$, we have that a basis of $H^2(\mathfrak{N}_2,{\bf k})$ is $\left\{[\Delta_{11}], [\Delta_{12}], [\Delta_{21}], [\Delta_{22}] \right\}$. Also, for $a, b, c, d\in {\bf k}$ we have
\begin{itemize}
\item $T_1(\mathfrak{N}_2)=\left\{a[\Delta_{11}]+[\Delta_{21}]+d[\Delta_{22}] \right\} \cup \left\{a[\Delta_{11}]+ [\Delta_{12}]+ c[\Delta_{21}]+d[\Delta_{22}]: c\neq 1 \right\}\cup \left\{[\Delta_{11}]+d[\Delta_{22}]: d\neq 0 \right\}\cup $ $ \cup\left\{a[\Delta_{11}]+ [\Delta_{12}]+ [\Delta_{21}]+d[\Delta_{22}]: ad\neq1 \right\}$.

\item $T_s(\mathfrak{N}_2)=Grass_s(H^2(\mathfrak{N}_2, {\bf k}))$ for $s=2, 3, 4$.
\end{itemize}

Additionally, the action of ${\mathcal A}(\mathfrak{N}_2,g)=GL(2,{\bf k})$ on
$\langle a[\Delta_{11}]+ b[\Delta_{12}]+c [\Delta_{21}]+ d[\Delta_{22}]\rangle \in  T_1(\mathfrak{N}_2)$ is the following:
$$\begin{array}{rclll}
(\alpha_{11}(a \alpha_{11}+ c \alpha_{21})+\alpha_{21}(b \alpha_{11}+ d \alpha_{21})) [\Delta_{11}]+
(\alpha_{12}(a \alpha_{11}+c\alpha_{21})+\alpha_{22}(b\alpha_{11}+d\alpha_{21}))[\Delta_{12}]+ \\
+(\alpha_{11}(a\alpha_{12}+c\alpha_{22})+\alpha_{21}(b\alpha_{12}+d\alpha_{22}))[\Delta_{21}]
+(\alpha_{12}(a\alpha_{12}+c\alpha_{22})+\alpha_{22}(b\alpha_{12}+d\alpha_{22}))[\Delta_{22}]
\end{array}$$

Again, we write each set $T_s(\mathfrak{N}_2)$ as a joint of their orbits, the following equalities can be verified by straightforward calculations:
\begin{enumerate}
    \item[$\star$] $T_1(\mathfrak{N}_2)=O( \langle [\Delta_{11}]+[\Delta_{22}]\rangle ) \cup O( \langle [\Delta_{12}]-[\Delta_{21}]\rangle ) \cup  O(\{ \langle \lambda [\Delta_{11}]+[\Delta_{21}]+[\Delta_{22}]\rangle:\lambda\in {\bf k}\}).$
        \item[$\star$] $T_2(\mathfrak{N}_2)= O( \langle [\Delta_{12}], [\Delta_{21}] \rangle ) \cup O( \langle [\Delta_{11}], [\Delta_{21}] \rangle ) \cup O(\langle [\Delta_{11}],[\Delta_{22}]\rangle)\cup O( \langle [\Delta_{11}]+[\Delta_{12}]+[\Delta_{22}], [\Delta_{21}]\rangle)\cup$

        $\cup O(\{\langle [\Delta_{11}]+[\Delta_{22}], [\Delta_{21}]+\lambda[\Delta_{22}]\rangle: \lambda\in {\bf k}_{\geq 0}\})\cup O(\{\langle [\Delta_{11}], [\Delta_{12}]+\lambda[\Delta_{21}]\rangle: \lambda\in {\bf k}\}).$
        \item[$\star$] $T_3(\mathfrak{N}_2)=O( \langle [\Delta_{12}], [\Delta_{21}], [\Delta_{22}]\rangle ) \cup  O( \langle [\Delta_{12}], [\Delta_{21}], [\Delta_{11}]+[\Delta_{22}]\rangle ) \cup O( \langle [\Delta_{11}], [\Delta_{12}]+[\Delta_{21}], [\Delta_{22}]\rangle ) \cup$

        $ \cup O(\{ \langle [\Delta_{11}]+[\Delta_{12}], [\Delta_{21}], \lambda [\Delta_{11}]+[\Delta_{22}]\rangle:\lambda \in {\bf k} \}).$
    \item[$\star$] $T_4(\mathfrak{N}_2)=H^2(\mathfrak{N}_2, {\bf k})$.
\end{enumerate}

Finally, we compute each radical extension of $(\mathfrak{N}_2, g)$ using its
correspondence with an orbit in $T_s(\mathfrak{N}_2):$
$$\begin{array}{ll llll}
(\mathfrak{N}_2)_{3,1} &:& e_1 e_1 = e_3 & e_1 e_2=0 & e_2 e_1=0  & e_2 e_2=e_3 \\
(\mathfrak{N}_2)_{3,2} &:& e_1 e_1 = 0 & e_1 e_2=e_3 & e_2 e_1=-e_3  & e_2 e_2=0 \\
(\mathfrak{N}_2)_{3,3}(\lambda \in {\bf k}) &:& e_1 e_1 = \lambda e_3 & e_1 e_2=0 & e_2 e_1=e_3  & e_2 e_2=e_3 \\
(\mathfrak{N}_2)_{4,4} &:& e_1 e_1 = 0 & e_1 e_2=e_3 & e_2 e_1=e_4  & e_2 e_2=0 \\
(\mathfrak{N}_2)_{4,5} &:& e_1 e_1 = e_3 & e_1 e_2=0 & e_2 e_1=e_4  & e_2 e_2=0 \\
(\mathfrak{N}_2)_{4,6} &:& e_1 e_1 = e_3 & e_1 e_2=0 & e_2 e_1=0  & e_2 e_2=e_4 \\
(\mathfrak{N}_2)_{4,7} &:& e_1 e_1 = e_3 & e_1 e_2=e_3 & e_2 e_1=e_4  & e_2 e_2=e_3 \\
(\mathfrak{N}_2)_{4,8}(\lambda \in {\bf k}_{\geq 0}) &:& e_1 e_1 = e_3 & e_1 e_2=0 & e_2 e_1=e_4  & e_2 e_2=e_3+\lambda e_4 \\
(\mathfrak{N}_2)_{4,9}(\lambda \in {\bf k}) &:& e_1 e_1 = e_3 & e_1 e_2=e_4 & e_2 e_1=\lambda e_4  & e_2 e_2=0 \\
(\mathfrak{N}_2)_{5,10} &:& e_1 e_1 = 0 & e_1 e_2=e_3 & e_2 e_1=e_4  & e_2 e_2=e_5 \\
(\mathfrak{N}_2)_{5,11} &:& e_1 e_1 = e_5 & e_1 e_2=e_3 & e_2 e_1=e_4  & e_2 e_2=e_5 \\
(\mathfrak{N}_2)_{5,12} &:& e_1 e_1 = e_3 & e_1 e_2=e_4 & e_2 e_1=e_4  & e_2 e_2=e_5 \\
(\mathfrak{N}_2)_{5,13} (\lambda \in {\bf k}) &:& e_1 e_1 = e_3+\lambda e_5 & e_1 e_2=e_3 & e_2 e_1=e_4  & e_2 e_2=e_5 \\
(\mathfrak{N}_2)_{6,14} &:& e_1 e_1 = e_3 & e_1 e_2=e_4 & e_2 e_1=e_5  & e_2 e_2=e_6
\end{array}$$

\newpage
\section{Appendix: Tables}
\begin{center}
\label{tab1}
{\bf Table 1.}  
\end{center}
\bigskip
{\tiny $$\begin{array}{lllll lll}
{\bf A}_{01}(\alpha)&:&({\bf A}_1 (\alpha\neq0))_{3,1} &:& 
e_1 e_1 = e_1+e_2 & e_1 e_2=\alpha e_2 & e_2 e_1=(1-\alpha)e_2+ e_3 & e_2 e_2=0 \\ 
{\bf A}_{02}&:&({\bf A}_1 (0))_{3,1} &:& 
e_1 e_1 = e_1+e_2 & e_1 e_2=e_3 & e_2 e_1=e_2 & e_2 e_2=0 \\ 
{\bf A}_{03}(\alpha)&:&({\bf A}_1 (\alpha\neq\frac{1}{2}))_{3,2} &:& 
e_1 e_1 = e_1+e_2 & e_1 e_2=\alpha e_2 & e_2 e_1=(1-\alpha)e_2 & e_2 e_2=e_3 \\
{\bf A}_{04}(\lambda)&:&({\bf A}_1 (\frac{1}{2}))_{3,2}(\lambda\in {\bf k}) &:& 
e_1 e_1 = e_1+e_2 & e_1 e_2=\frac{1}{2} e_2 & e_2 e_1=\frac{1}{2} e_2 +\lambda e_3 & e_2 e_2=e_3 \\

{\bf A}_{05}&:&({\bf A}_2 )_{3,1} &:&    e_1 e_1 = e_2+e_3 & e_1 e_2= e_2 & e_2 e_1=-e_2  & e_2 e_2=0\\
{\bf A}_{06}&:&({\bf A}_2 )_{3,2} &:&    e_1 e_1 = e_2 & e_1 e_2= e_2 & e_2 e_1=-e_2+e_3  & e_2 e_2=0\\
{\bf A}_{07}(\lambda)&:&({\bf A}_2 )_{3,3}(\lambda\in {\bf k}) &:&    e_1 e_1 = e_2+\lambda e_3 & e_1 e_2= e_2 & e_2 e_1=-e_2  & e_2 e_2=e_3 \\
{\bf A}_{08}&:&({\bf A}_3 )_{3,1} &:& e_1 e_1 = e_2 & e_1 e_2=0 & e_2 e_1=0 & e_2 e_2=e_3 \\
{\bf A}_{09}&:&({\bf A}_3 )_{3,2} &:& e_1 e_1 = e_2 & e_1 e_2=0 & e_2 e_1= e_3 & e_2 e_2=e_3 \\
{\bf A}_{10}&:&({\bf A}_3 )_{3,3} &:& e_1 e_1 = e_2 & e_1 e_2=0 & e_2 e_1=e_3 & e_2 e_2=0 \\
{\bf A}_{11}(\lambda)&:&({\bf A}_3 )_{3,4}(\lambda\in {\bf k}) &:& e_1 e_1 = e_2 & e_1 e_2=e_3 & e_2 e_1=\lambda e_3 & e_2 e_2=0 \\

{\bf A}_{12}(\alpha)&:&({\bf A}_4 (\alpha))_{3,1} &:& e_1 e_1 = \alpha e_1+e_2 & e_1 e_2= e_1+\alpha e_2 & e_2 e_1=-e_1  & e_2 e_2= e_3 \\
{\bf A}_{13}(\alpha, \lambda)&:&({\bf A}_4 (\alpha\neq 0))_{3,2}(\lambda\in{\bf k}) &:& e_1 e_1 = \alpha e_1+e_2 & e_1 e_2= e_1+\alpha e_2 + e_3 & e_2 e_1=-e_1   & e_2 e_2=\lambda e_3 \\
{\bf A}_{14}(\lambda)&:&({\bf A}_4 (0))_{3,2}(\lambda\in{\bf k}_{\geq 0}\}) &:& e_1 e_1 = e_2 & e_1 e_2= e_1 + e_3 & e_2 e_1=-e_1   & e_2 e_2=\lambda e_3 \\

{\bf A}_{15}(\alpha)&:&({\bf B}_1(\alpha))_{3,1} &:& e_1 e_1 = 0 & e_1 e_2=(1- \alpha) e_1 + e_2  & e_2 e_1=\alpha e_1- e_2 &  e_2 e_2= e_3  \\
{\bf A}_{16}(\alpha, \lambda)&:&({\bf B}_1(\alpha))_{3,2}(\lambda\in {\bf k}) &:& e_1 e_1 = e_3 & e_1 e_2=(1- \alpha) e_1 + e_2  & e_2 e_1=\alpha e_1- e_2 &  e_2 e_2= \lambda  e_3  \\
{\bf A}_{17}(\alpha)&:&({\bf B}_2(\alpha))_{3,1} &:&
e_1 e_1 = 0  & e_1 e_2=(1- \alpha) e_1 & e_2 e_1=\alpha e_1 &  e_2 e_2=e_3 \\
{\bf A}_{18}(\alpha)&:&({\bf B}_2(\alpha\neq0))_{3,2} &:&
e_1 e_1 = 0  & e_1 e_2=(1- \alpha) e_1+e_3 & e_2 e_1=\alpha e_1 &  e_2 e_2=e_3 \\
{\bf A}_{19}&:&({\bf B}_2(0))_{3,2} &:&
e_1 e_1 = 0  & e_1 e_2=e_1 & e_2 e_1=e_3 &  e_2 e_2=e_3 \\
{\bf A}_{20}(\alpha)&:&({\bf B}_2(\alpha\neq0))_{3,3} &:&
e_1 e_1 = 0  & e_1 e_2=(1- \alpha) e_1+e_3 & e_2 e_1=\alpha e_1 &  e_2 e_2=0 \\
{\bf A}_{21}&:&({\bf B}_2(0))_{3,3} &:&
e_1 e_1 = 0  & e_1 e_2= e_1 & e_2 e_1=e_3 &  e_2 e_2=0 \\
{\bf A}_{22}(\alpha)&:&({\bf B}_2(\alpha))_{3,4} &:&
e_1 e_1 = e_3  & e_1 e_2=(1- \alpha) e_1 & e_2 e_1=\alpha e_1 &  e_2 e_2=0 \\
{\bf A}_{23}(\alpha)&:&({\bf B}_2(\alpha))_{3,5} &:&
e_1 e_1 = e_3  & e_1 e_2=(1- \alpha) e_1 & e_2 e_1=\alpha e_1 &  e_2 e_2=e_3 \\
{\bf A}_{24}(\alpha, \lambda)&:&({\bf B}_2(\alpha\neq0))_{3,6}(\lambda\in {\bf k}) &:&
e_1 e_1 = e_3  & e_1 e_2=(1- \alpha) e_1+e_3 & e_2 e_1=\alpha e_1 &  e_2 e_2=\lambda e_3 \\
{\bf A}_{25}(\lambda)&:&({\bf B}_2(0))_{3,6}(\lambda\in {\bf k}) &:&
e_1 e_1 = e_3  & e_1 e_2=e_1 & e_2 e_1=e_3 &  e_2 e_2=\lambda e_3 \\

{\bf A}_{26}&:&({\bf B}_3)_{3,1} &:& e_1 e_1 = e_3  & e_1 e_2= e_2 & e_2 e_1=-e_2 & e_2 e_2=0\\
{\bf A}_{27}&:&({\bf B}_3)_{3,2} &:& e_1 e_1 = 0  & e_1 e_2= e_2 & e_2 e_1=-e_2+e_3 & e_2 e_2=0\\
{\bf A}_{28}&:&({\bf B}_3)_{3,3} &:& e_1 e_1 = e_3  & e_1 e_2= e_2 & e_2 e_1=-e_2 & e_2 e_2=e_3\\
{\bf A}_{29}&:&({\bf B}_3)_{3,4} &:& e_1 e_1 = 0  & e_1 e_2= e_2 & e_2 e_1=-e_2 & e_2 e_2=e_3\\

{\bf A}_{30}(\alpha, \beta)&:&({\bf C}(\alpha, \beta))_{3,1} &:&  e_1 e_1 = e_2 & e_1 e_2=(1-\alpha)e_1+ \beta e_2  & e_2 e_1=\alpha e_1 - \beta e_2 & e_2 e_2=e_2+e_3 \\

{\bf A}_{31}(\alpha, \lambda)&:&({\bf C}(\alpha\neq 0, 0))_{3,2}\\
&&(\lambda\in{\bf k}_{\geq 0}) &:&  e_1 e_1 = e_2 & e_1 e_2=(1-\alpha)e_1 + e_3  & e_2 e_1=\alpha e_1 & e_2 e_2=e_2+\lambda e_3 \\

{\bf A}_{32}(\lambda)&:&({\bf C}(0, 0))_{3,2}(\lambda\in{\bf k}_{\geq 0}) &:&  e_1 e_1 = e_2 & e_1 e_2= e_1 & e_2 e_1= e_3  & e_2 e_2=e_2+\lambda e_3 \\

{\bf A}_{33}(\alpha, \beta, \lambda)&:&({\bf C}(\alpha\neq 0, \beta\neq0))_{3,2}\\
&&(\lambda\in{\bf k}) &:&  e_1 e_1 = e_2 & e_1 e_2=(1-\alpha)e_1+ \beta e_2+e_3 & e_2 e_1=\alpha e_1 - \beta e_2 & e_2 e_2=e_2+\lambda e_3 \\

{\bf A}_{34}(\beta, \lambda)&:&({\bf C}(0, \beta\neq0))_{3,2}(\lambda\in{\bf k}) &:&  e_1 e_1 = e_2 & e_1 e_2=  e_1 + \beta e_2 & e_2 e_1= -\beta e_2+e_3 & e_2 e_2=e_2+\lambda e_3 \\

{\bf A}_{35}(\alpha, \beta)&:&({\bf D}_1(\alpha,\beta))_{3,1}  &:& 
e_1 e_1 = e_1 & e_1 e_2=(1-\alpha)e_1+ \beta e_2  & e_2 e_1=\alpha e_1 - \beta e_2 & e_2 e_2=e_3\\
{\bf A}_{36}(\alpha, \lambda)&:&({\bf D}_1(\alpha\neq\frac{1}{2},0))_{3,2}(\lambda \in {\bf k}) &:& 
e_1 e_1 = e_1 & e_1 e_2=(1-\alpha)e_1+e_3  & e_2 e_1=\alpha e_1 &e_2 e_2=\lambda e_3 \\
{\bf A}_{37}(\lambda)&:&({\bf D}_1(\frac{1}{2},0)_{3,2}((\lambda, 0) \in \mathcal U) &:& 
e_1 e_1 = e_1 & e_1 e_2=\frac{1}{2}e_1+e_3  & e_2 e_1=\frac{1}{2} e_1 &e_2 e_2=\lambda e_3 \\
{\bf A}_{38}(\alpha, \lambda, \mu)&:&({\bf D}_1(\alpha\neq\frac{1}{2},0))_{3,3}\\
&&(\lambda, \mu \in {\bf k}) &:& 
e_1 e_1 = e_1 & e_1 e_2=(1-\alpha)e_1+\lambda e_3  & e_2 e_1=\alpha e_1+e_3 &e_2 e_2=\mu e_3 \\
{\bf A}_{39}(\mu, \lambda)&:&({\bf D}_1(\frac{1}{2},0))_{3,3}\\
&&((\mu, \lambda) \in \mathcal U) &:& 
e_1 e_1 = e_1 & e_1 e_2=\frac{1}{2}e_1+\lambda e_3  & e_2 e_1=\frac{1}{2} e_1+e_3 &e_2 e_2=\mu e_3 \\

{\bf A}_{40}(\alpha,\lambda)&:&({\bf D}_1(\alpha,2\alpha-1\neq0))_{3,4}\\

&& ((\lambda, 0)\in \mathcal U)&:& e_1 e_1 = e_1 & e_1 e_2=(1-\alpha)e_1+ (2\alpha-1) e_2  & e_2 e_1=\alpha e_1 - (2\alpha-1) e_2+ e_3 & e_2 e_2=\lambda e_3\\

{\bf A}_{41}(\alpha, \beta, \lambda)&:&({\bf D}_1(\alpha,\beta))_{3,4}(\lambda \in {\bf k})\\

&&(\beta\neq2\alpha-1, 0) &:& 
e_1 e_1 = e_1 & e_1 e_2=(1-\alpha)e_1+ \beta e_2  & e_2 e_1=\alpha e_1 - \beta e_2+ e_3 & e_2 e_2=\lambda e_3\\

{\bf A}_{42}(\alpha, \beta)&:&({\bf D}_2(\alpha,\beta))_{3,1} &:& e_1 e_1 = e_1 & e_1 e_2=\alpha e_2  & e_2 e_1= \beta e_2 & e_2 e_2=e_3 \\
{\bf A}_{43}(\alpha)&:&({\bf D}_2(\alpha,0)_{3,2} &:& e_1 e_1 = e_1 & e_1 e_2=\alpha e_2 & e_2 e_1= e_3 & e_2 e_2=e_3 \\

{\bf A}_{44}(\alpha)&:&({\bf D}_2(\alpha, 0)_{3,3} &:& e_1 e_1 = e_1 & e_1 e_2=\alpha e_2 & e_2 e_1=e_3   & e_2 e_2=0 \\

{\bf A}_{45}( \lambda)&:&({\bf D}_2(0,0))_{3,4}(\lambda\in {\bf k}) &:& e_1 e_1 = e_1 & e_1 e_2=e_3 & e_2 e_1=\lambda e_3 & e_2 e_2=e_3\\
{\bf A}_{46}( \lambda)&:&({\bf D}_2(0,0))_{3,5}(\lambda\in {\bf k}) &:& e_1 e_1 = e_1 & e_1 e_2=e_3 & e_2 e_1=\lambda e_3 & e_2 e_2=0\\
{\bf A}_{47}(\alpha, \beta)&:&({\bf D}_2(\alpha,\beta\neq0)_{3,6} &:& e_1 e_1 = e_1 & e_1 e_2=\alpha e_2 +e_3 & e_2 e_1= \beta e_2 & e_2 e_2=e_3 \\

{\bf A}_{48}(\alpha, \beta)&:&({\bf D}_2(\alpha,\beta\neq0)_{3,7} &:& e_1 e_1 = e_1 & e_1 e_2=\alpha e_2+e_3 & e_2 e_1= \beta e_2  & e_2 e_2=0 \\

{\bf A}_{49}(\alpha, \beta)&:&({\bf D}_3(\alpha,\beta))_{3,1} &:& 
e_1 e_1 = e_1 & e_1 e_2=e_1+ \alpha e_2   & e_2 e_1=- e_1 + \beta e_2 &  e_2 e_2=e_3 \\ 
{\bf A}_{50}(\beta, \lambda)&:&({\bf D}_3(0,\beta))_{3,2}(\lambda\in {\bf k}) &:& 
e_1 e_1 = e_1 & e_1 e_2= e_1 +e_3   & e_2 e_1= -e_1+\beta e_2   &  e_2 e_2=\lambda e_3 \\ 
{\bf A}_{51}(\lambda, \mu)&:&({\bf D}_3(0,0))_{3,3}(\lambda, \mu\in {\bf k}) &:& 
e_1 e_1 = e_1 & e_1 e_2=e_1+\lambda e_3    & e_2 e_1=- e_1+e_3 &  e_2 e_2=\mu e_3 \\ 
{\bf A}_{52}(\alpha, \beta, \lambda)&:&({\bf D}_3(\alpha\neq0,\beta))_{3,4}\\
&&(\lambda\in {\bf k}) &:& 
e_1 e_1 = e_1 & e_1 e_2=e_1+ \alpha e_2   & e_2 e_1=- e_1 + \beta e_2+ e_3 &  e_2 e_2=\lambda e_3 \\

{\bf A}_{53}(\alpha, \beta, \gamma,\delta)&:&({\bf E}_1(\alpha, \beta, \gamma,\delta))_{3,1} &:&  e_1 e_1 = e_1 & e_1 e_2=\alpha e_1+ \beta e_2 & e_2 e_1= \gamma  e_1 + \delta  e_2 +e_3& e_2 e_2=e_2 \\

{\bf A}_{54}(\alpha, \beta, \gamma,\delta, \lambda)&:&({\bf E}_1(\alpha, \beta, \gamma,\delta))_{3,2}\\
&&(\lambda \in {\bf k})\\
&&(\gamma\neq\beta \textrm{ or } \delta\neq\alpha) &:&  e_1 e_1 = e_1 & e_1 e_2=\alpha e_1+ \beta e_2+e_3 & e_2 e_1= \gamma  e_1 + \delta  e_2 +\lambda e_3& e_2 e_2=e_2 \\

{\bf A}_{55}(\alpha, \beta, \lambda)&:&({\bf E}_1(\alpha, \beta, \beta,\alpha))_{3,2}\\
&&(\lambda \in {\bf k}^*_{>1}\cup \{1\}) &:&   e_1 e_1 = e_1 & e_1 e_2=\alpha e_1+ \beta e_2+e_3 & e_2 e_1= \beta  e_1 + \alpha  e_2 +\lambda e_3& e_2 e_2=e_2 \\

{\bf A}_{56}(\alpha, \beta, \gamma)&:&({\bf E}_2(\alpha, \beta, \gamma))_{3,1} &:&  e_1 e_1 = e_1 & e_1 e_2=(1-\alpha)e_1+ \beta e_2 & e_2 e_1= \alpha e_1 + \gamma e_2 +e_3& e_2 e_2=e_2 \\
{\bf A}_{57}(\alpha, \beta, \gamma, \lambda)&:&({\bf E}_2(\alpha, \beta, \gamma))_{3,2}(\lambda\in {\bf k}) &:&  e_1 e_1 = e_1 & e_1 e_2=(1-\alpha)e_1+ \beta e_2+e_3 & e_2 e_1= \alpha e_1 + \gamma e_2 +\lambda e_3& e_2 e_2=e_2 \\

{\bf A}_{58}(\alpha, \beta, \gamma)&:&({\bf E}_3(\alpha, \beta, \gamma))_{3,1} &:&  e_1 e_1 = e_1 & e_1 e_2=(1-\alpha)\gamma e_1+ \frac{\beta}{\gamma} e_2  & e_2 e_1=\alpha \gamma e_1 + \frac{1-\beta}{\gamma}e_2 +e_3&
e_2 e_2=e_2 \\
{\bf A}_{59}(\alpha, \beta, \gamma, \lambda)&:&({\bf E}_3(\alpha, \beta, \gamma))_{3,2}(\lambda\in {\bf k}) \\
&&(\beta\neq \alpha \textrm{ or } \gamma\neq -1) &:&  e_1 e_1 = e_1 & e_1 e_2=(1-\alpha)\gamma e_1+ \frac{\beta}{\gamma} e_2  +e_3& e_2 e_1=\alpha \gamma e_1 + \frac{1-\beta}{\gamma}e_2+\lambda e_3 &
e_2 e_2=e_2\\
{\bf A}_{60}(\alpha,\lambda)&:&({\bf E}_3(\alpha, \alpha, -1))_{3,2} &&&\\

&&(\lambda\in {\bf k}^*_{>1}\cup \{1\}) &:&  e_1 e_1 = e_1 & e_1 e_2=-(1-\alpha) e_1-\alpha e_2  +e_3& e_2 e_1=-\alpha e_1  -(1-\alpha)e_2+\lambda e_3 &
e_2 e_2=e_2 \\

{\bf A}_{61}&:&({\bf E}_4)_{3,1} &:&  e_1 e_1 = e_1 & e_1 e_2=e_1+ e_2  &  e_2 e_1=e_3  & e_2 e_2=e_2 \\
{\bf A}_{62}(\lambda)&:&({\bf E}_4)_{3,2}(\lambda\in {\bf k}) &:&  e_1 e_1 = e_1 & e_1 e_2=e_1+ e_2+e_3  &  e_2 e_1=\lambda e_3  & e_2 e_2=e_2 \\

{\bf A}_{63}(\alpha)&:&({\bf E}_5(\alpha))_{3,1} &:& e_1 e_1 = e_1 & e_1 e_2=(1-\alpha)e_1+ \alpha e_2 & e_2 e_1=\alpha e_1+ (1-\alpha)e_2+e_3   & e_2 e_2=e_2 \\
{\bf A}_{64}(\alpha)&:&({\bf E}_5(\alpha))_{3,2} &:& e_1 e_1 = e_1 & e_1 e_2=(1-\alpha)e_1+ \alpha e_2+e_3 & e_2 e_1=\alpha e_1+ (1-\alpha) e_2-e_3  & e_2 e_2=e_2 \\
{\bf A}_{65}&:&({\bf E}_5(1/2))_{3,3} &:& e_1 e_1 = e_1 & e_1 e_2=\frac{1}{2}e_1+ \frac{1}{2} e_2+e_3 & e_2 e_1=\frac{1}{2} e_1+ \frac{1}{2} e_2+e_3  & e_2 e_2=e_2 \\

{\bf A}_{66}&:&(\mathfrak{N}_2)_{3,1} &:& e_1 e_1 = e_3 & e_1 e_2=0 & e_2 e_1=0  & e_2 e_2=e_3 \\

{\bf A}_{67}&:&(\mathfrak{N}_2)_{3,2} &:& e_1 e_1 = 0 & e_1 e_2=e_3 & e_2 e_1=-e_3  & e_2 e_2=0 \\

{\bf A}_{68}(\lambda)&:&(\mathfrak{N}_2)_{3,3}(\lambda \in {\bf k}) &:& e_1 e_1 = \lambda e_3 & e_1 e_2=0 & e_2 e_1=e_3  & e_2 e_2=e_3 
\end{array}$$} 

\newpage

\begin{center}
\label{tab2}
{\bf Table 2.} 
\end{center}
\bigskip
{\tiny $$\begin{array}{lllll lll}

{\bf A}_{69}(\alpha) &:&({\bf A}_1 (\alpha\neq0))_{4,3} &:& 
e_1 e_1 = e_1+e_2 & e_1 e_2=\alpha e_2 & e_2 e_1=(1-\alpha)e_2+e_3 & e_2 e_2=e_4 \\

{\bf A}_{70} &:&({\bf A}_1 (0))_{4,3} &:& 
e_1 e_1 = e_1+e_2 & e_1 e_2= e_3 & e_2 e_1=e_2 & e_2 e_2=e_4 \\

{\bf A}_{71} &:&({\bf A}_2 )_{4,4} &:&    e_1 e_1 = e_2 & e_1 e_2= e_2 & e_2 e_1=-e_2+e_3  & e_2 e_2=e_4\\

{\bf A}_{72}(\lambda) &:&({\bf A}_2 )_{4,5}(\lambda\in {\bf k}) &:&    e_1 e_1 = e_2+e_3 & e_1 e_2= e_2 & e_2 e_1=-e_2+e_4  & e_2 e_2=\lambda e_3\\

{\bf A}_{73} &:&({\bf A}_2 )_{4,6} &:&    e_1 e_1 = e_2+e_3 & e_1 e_2= e_2 & e_2 e_1=-e_2 & e_2 e_2=e_4\\

{\bf A}_{74} &:&({\bf A}_3 )_{4,5} &:& e_1 e_1 = e_2 & e_1 e_2=0 & e_2 e_1=e_4 & e_2 e_2=e_3 \\
{\bf A}_{75} &:&({\bf A}_3 )_{4,6} &:& e_1 e_1 = e_2 & e_1 e_2=e_4 & e_2 e_1=e_3 & e_2 e_2=0 \\
{\bf A}_{76}(\lambda) &:&({\bf A}_3 )_{4,7} (\lambda\in {\bf k}) &:& e_1 e_1 = e_2 & e_1 e_2=e_4 & e_2 e_1=\lambda e_4 & e_2 e_2=e_3 \\

{\bf A}_{77}(\alpha) &:&({\bf A}_4 (\alpha))_{4,3}&:& e_1 e_1 = \alpha e_1+e_2 & e_1 e_2= e_1+\alpha e_2+ e_3 & e_2 e_1=-e_1  & e_2 e_2= e_4 \\

{\bf A}_{78}(\alpha) &:&({\bf B}_1(\alpha))_{4,3} &:& e_1 e_1 = e_3 & e_1 e_2=(1- \alpha) e_1 + e_2  & e_2 e_1=\alpha e_1- e_2 &  e_2 e_2= e_4  \\

{\bf A}_{79}(\alpha) &:&({\bf B}_2(\alpha\neq0))_{4,7} &:&
e_1 e_1 = 0  & e_1 e_2=(1- \alpha) e_1+e_3 & e_2 e_1=\alpha e_1 &  e_2 e_2=e_4 \\
{\bf A}_{80} &:&({\bf B}_2(0))_{4,7} &:&
e_1 e_1 = 0  & e_1 e_2= e_1 & e_2 e_1=e_3 &  e_2 e_2=e_4 \\
{\bf A}_{81}(\alpha) &:&({\bf B}_2(\alpha\neq0))_{4,8} &:&
e_1 e_1 = e_3  & e_1 e_2=(1- \alpha) e_1+e_4 & e_2 e_1=\alpha e_1 &  e_2 e_2=e_3 \\
{\bf A}_{82} &:&({\bf B}_2(0))_{4,8} &:&
e_1 e_1 = e_3  & e_1 e_2=e_1 & e_2 e_1=e_4 &  e_2 e_2=e_3 \\
{\bf A}_{83}(\alpha) &:&({\bf B}_2(\alpha\neq0))_{4,9} &:&
e_1 e_1 = e_3  & e_1 e_2=(1- \alpha) e_1+e_4 & e_2 e_1=\alpha e_1 &  e_2 e_2=0 \\
{\bf A}_{84} &:&({\bf B}_2(0))_{4,9} &:&
e_1 e_1 = e_3  & e_1 e_2= e_1 & e_2 e_1=e_4 &  e_2 e_2=0 \\
{\bf A}_{85}(\alpha) &:&({\bf B}_2(\alpha\neq0))_{4,10} &:&
e_1 e_1 = e_3  & e_1 e_2=(1- \alpha) e_1+e_3 & e_2 e_1=\alpha e_1 &  e_2 e_2=e_4 \\
{\bf A}_{86} &:&({\bf B}_2(0))_{4,10} &:&
e_1 e_1 = e_3  & e_1 e_2=e_1 & e_2 e_1=e_3 &  e_2 e_2=e_4 \\
{\bf A}_{87}(\alpha) &:&({\bf B}_2(\alpha))_{4,11} &:&
e_1 e_1 = e_3  & e_1 e_2=(1- \alpha) e_1 & e_2 e_1=\alpha e_1 &  e_2 e_2=e_4 \\
{\bf A}_{88}(\alpha, \lambda) &:&({\bf B}_2(\alpha\neq0))_{4,12}(\lambda\in {\bf k}) &:&
e_1 e_1 = e_3  & e_1 e_2=(1- \alpha) e_1+ \lambda e_3+e_4 & e_2 e_1=\alpha e_1 &  e_2 e_2=e_4 \\
{\bf A}_{89}(\lambda) &:&({\bf B}_2(0))_{4,12}(\lambda\in {\bf k}) &:&
e_1 e_1 = e_3  & e_1 e_2=e_1 & e_2 e_1=\lambda  e_3+e_4 &  e_2 e_2=e_4 \\

{\bf A}_{90} &:&({\bf B}_3)_{4,5} &:& e_1 e_1 = 0  & e_1 e_2= e_2 & e_2 e_1=-e_2+e_3 & e_2 e_2=e_4\\
{\bf A}_{91} &:&({\bf B}_3)_{4,6} &:& e_1 e_1 = e_3  & e_1 e_2= e_2 & e_2 e_1=-e_2+e_4 & e_2 e_2=e_3\\
{\bf A}_{92} &:&({\bf B}_3)_{4,7} &:& e_1 e_1 = e_3  & e_1 e_2= e_2 & e_2 e_1=-e_2+e_4 & e_2 e_2=0\\
{\bf A}_{93} &:&({\bf B}_3)_{4,8} &:& e_1 e_1 = e_3  & e_1 e_2= e_2 & e_2 e_1=-e_2 & e_2 e_2=e_4\\

{\bf A}_{94}(\alpha, \beta, \lambda) &:&({\bf C}(\alpha\neq 0, \beta))_{4,3}(\lambda\in{\bf k}) &:&  e_1 e_1 = e_2 & e_1 e_2=(1-\alpha)e_1+ \beta e_2+e_3  & e_2 e_1=\alpha e_1 - \beta e_2  & e_2 e_2=e_2+e_4 \\
{\bf A}_{95}(\beta, \lambda) &:&({\bf C}(0, \beta))_{4,3}(\lambda\in{\bf k}) &:&  e_1 e_1 = e_2 & e_1 e_2= e_1 + \beta e_2  & e_2 e_1=- \beta e_2 +e_3  & e_2 e_2=e_2+e_4 \\

{\bf A}_{96}(\alpha, \beta) &:&({\bf D}_1(\alpha,\beta))_{4,5}  &:& 
e_1 e_1 = e_1 & e_1 e_2=(1-\alpha)e_1+ \beta e_2  & e_2 e_1=\alpha e_1 - \beta e_2+e_3 & e_2 e_2=e_4  \\
{\bf A}_{97}(\alpha, \lambda) &:&({\bf D}_1(\alpha,0))_{4,6}(\lambda \in {\bf k}) &:& 
e_1 e_1 = e_1 & e_1 e_2=(1-\alpha)e_1 +e_4 & e_2 e_1=\alpha e_1+\lambda e_4 &e_2 e_2= e_3 \\

{\bf A}_{98}(\alpha, \lambda, \mu) &:&({\bf D}_1(\alpha\neq\frac{1}{2},0))_{4,7}\\

&&(\lambda,\mu \in {\bf k}) &:& 
e_1 e_1 = e_1 & e_1 e_2=(1-\alpha)e_1 +e_4 & e_2 e_1=\alpha e_1+ e_3 &e_2 e_2= \lambda e_3+\mu e_4 \\

{\bf A}_{99}(\lambda, \mu) &:&({\bf D}_1(\frac{1}{2},0))_{4,7}\\
&&((\lambda, 0), (\mu,0) \in \mathcal U) &:& 
e_1 e_1 = e_1 & e_1 e_2=\frac{1}{2}e_1 +e_4 & e_2 e_1=\frac{1}{2} e_1+ e_3 &e_2 e_2= \lambda e_3+\mu e_4 \\

{\bf A}_{100}(\alpha) &:&({\bf D}_2(\alpha, 0)_{4,8} &:& e_1 e_1 = e_1 & e_1 e_2=\alpha e_2 & e_2 e_1= e_3  & e_2 e_2=e_4 \\

{\bf A}_{101} &:&({\bf D}_2(0,0))_{4,9} &:& e_1 e_1 = e_1 & e_1 e_2=e_3 & e_2 e_1= e_4 & e_2 e_2=e_4 \\
{\bf A}_{102} &:&({\bf D}_2(0,0))_{4,10} &:& e_1 e_1 = e_1 & e_1 e_2=e_3 & e_2 e_1= e_4 & e_2 e_2=0 \\
{\bf A}_{103}( \lambda) &:&({\bf D}_2(0,0))_{4,11}(\lambda\in {\bf k}) &:& e_1 e_1 = e_1 & e_1 e_2=e_4 & e_2 e_1= \lambda e_4 & e_2 e_2=e_3 \\
{\bf A}_{104}( \lambda) &:&({\bf D}_2(0,0))_{4,12}(\lambda\in {\bf k}) &:& e_1 e_1 = e_1 & e_1 e_2=e_4 & e_2 e_1= e_3 & e_2 e_2=\lambda e_3+e_4 \\
{\bf A}_{105}(\alpha, \beta) &:&({\bf D}_2(\alpha,\beta\neq0)_{4,13} &:& e_1 e_1 = e_1 & e_1 e_2=\alpha e_2+e_3 & e_2 e_1= \beta e_2  & e_2 e_2=e_4 \\

{\bf A}_{106} &:&({\bf D}_3(0,0))_{4,5} &:& 
e_1 e_1 = e_1 & e_1 e_2=e_1   & e_2 e_1=- e_1+e_3 &  e_2 e_2=e_4 \\ 
{\bf A}_{107}(\alpha, \beta) &:&({\bf D}_3(\alpha\neq0,\beta))_{4,5} &:& 
e_1 e_1 = e_1 & e_1 e_2=e_1+ \alpha e_2   & e_2 e_1=- e_1 + \beta e_2+e_3 &  e_2 e_2=e_4 \\ 
{\bf A}_{108}(\lambda) &:&({\bf D}_3(0,0))_{4,6}(\lambda \in {\bf k}) &:& 
e_1 e_1 = e_1 & e_1 e_2=e_1+e_4    & e_2 e_1=- e_1+\lambda e_4 &  e_2 e_2=e_3 \\ 
{\bf A}_{109}(\lambda, \mu) &:&({\bf D}_3(0,0))_{4,7}(\lambda, \mu\in {\bf k}) &:& 
e_1 e_1 = e_1 & e_1 e_2=e_1 + e_4    & e_2 e_1=- e_1+e_3 &  e_2 e_2=\lambda e_3 +\mu e_4 \\ 
{\bf A}_{110}(\beta) &:&({\bf D}_3(0,\beta\neq0))_{4,8} &:& 
e_1 e_1 = e_1 & e_1 e_2= e_1 +e_3     & e_2 e_1= -e_1+\beta e_2       &  e_2 e_2=e_4 \\ 

{\bf A}_{111}(\alpha, \beta, \gamma,\delta) &:&({\bf E}_1(\alpha, \beta, \gamma,\delta))_{4,3} &:&  e_1 e_1 = e_1 & e_1 e_2=\alpha e_1+ \beta e_2+e_3 & e_2 e_1= \gamma  e_1 + \delta  e_2 + e_4 & e_2 e_2=e_2 \\

{\bf A}_{112}(\alpha, \beta, \gamma) &:&({\bf E}_2(\alpha, \beta, \gamma))_{4,3} &:&  e_1 e_1 = e_1 & e_1 e_2=(1-\alpha)e_1+ \beta e_2+e_3 & e_2 e_1= \alpha e_1 + \gamma e_2+e_4 & e_2 e_2=e_2 \\

{\bf A}_{113}(\alpha, \beta, \gamma) &:&({\bf E}_3(\alpha, \beta, \gamma))_{4,3} &:&  e_1 e_1 = e_1 & e_1 e_2=(1-\alpha)\gamma e_1+ \frac{\beta}{\gamma} e_2+e_3  & e_2 e_1=\alpha \gamma e_1 + \frac{1-\beta}{\gamma}e_2 +e_4&
e_2 e_2=e_2 \\

{\bf A}_{114} &:&({\bf E}_4)_{4,3} &:&  e_1 e_1 = e_1 & e_1 e_2=e_1+ e_2+e_3  &  e_2 e_1= e_4  & e_2 e_2=e_2 \\

{\bf A}_{115}(\alpha) &:&({\bf E}_5(\alpha))_{4,4} &:& e_1 e_1 = e_1 & e_1 e_2=(1-\alpha)e_1+ \alpha e_2+e_3 & e_2 e_1=\alpha e_1+ (1-\alpha) e_2+e_4  & e_2 e_2=e_2 \\

{\bf A}_{116}&:&(\mathfrak{N}_2)_{4,4} &:& e_1 e_1 = 0 & e_1 e_2=e_3 & e_2 e_1=e_4  & e_2 e_2=0 \\
{\bf A}_{117}&:&(\mathfrak{N}_2)_{4,5} &:& e_1 e_1 = e_3 & e_1 e_2=0 & e_2 e_1=e_4  & e_2 e_2=0 \\
{\bf A}_{118}&:&(\mathfrak{N}_2)_{4,6} &:& e_1 e_1 = e_3 & e_1 e_2=0 & e_2 e_1=0  & e_2 e_2=e_4 \\
{\bf A}_{119}&:&(\mathfrak{N}_2)_{4,7} &:& e_1 e_1 = e_3 & e_1 e_2=e_3 & e_2 e_1=e_4  & e_2 e_2=e_3 \\
{\bf A}_{120}(\lambda)&:&(\mathfrak{N}_2)_{4,8}(\lambda \in {\bf k}_{\geq 0}) &:& e_1 e_1 = e_3 & e_1 e_2=0 & e_2 e_1=e_4  & e_2 e_2=e_3+\lambda e_4 \\
{\bf A}_{121}(\lambda)&:&(\mathfrak{N}_2)_{4,9}(\lambda \in {\bf k}) &:& e_1 e_1 = e_3 & e_1 e_2=e_4 & e_2 e_1=\lambda e_4  & e_2 e_2=0 \\

\end{array}$$}

\bigskip
\bigskip

\begin{center}
\label{tab3}
{\bf Table 3.} 
\end{center}
\bigskip
$$\begin{array}{lllll lll}
{\bf A}_{122}&:&({\bf A}_2 )_{5,7}  &:&    e_1 e_1 = e_2+ e_3 & e_1 e_2= e_2 & e_2 e_1=-e_2+e_4  & e_2 e_2=e_5 \\
{\bf A}_{123}&:&({\bf A}_3 )_{5,8} &:& e_1 e_1 = e_2 & e_1 e_2=e_3 & e_2 e_1=e_4 & e_2 e_2=e_5 \\

{\bf A}_{124}(\alpha)&:&({\bf B}_2(\alpha\neq0))_{5,13} &:&
e_1 e_1 = e_3  & e_1 e_2=(1- \alpha) e_1+e_4 & e_2 e_1=\alpha e_1&  e_2 e_2=e_5 \\
{\bf A}_{125}&:&({\bf B}_2(0))_{5,13} &:&
e_1 e_1 = e_3  & e_1 e_2= e_1 & e_2 e_1=e_4 &  e_2 e_2=e_5 \\

{\bf A}_{126}&:&({\bf B}_3)_{5,9} &:& e_1 e_1 = e_3  & e_1 e_2= e_2 & e_2 e_1=-e_2+e_4 & e_2 e_2=e_5\\

{\bf A}_{127}(\alpha)&:&({\bf D}_1(\alpha,0))_{5,9} &:& 
e_1 e_1 = e_1 & e_1 e_2=(1-\alpha)e_1+e_3  & e_2 e_1=\alpha e_1+e_4 &e_2 e_2=e_5 \\  

{\bf A}_{128}&:&({\bf D}_2(0, 0))_{5,14} &:& e_1 e_1 = e_1 & e_1 e_2=e_3 & e_2 e_1= e_4 & e_2 e_2=e_5\\

{\bf A}_{129}&:&({\bf D}_3(0,0))_{5,9} &:& 
e_1 e_1 = e_1 & e_1 e_2=e_1+e_3    & e_2 e_1=- e_1+e_4 &  e_2 e_2=e_5\\

{\bf A}_{130}&:&(\mathfrak{N}_2)_{5,10} &:& e_1 e_1 = 0 & e_1 e_2=e_3 & e_2 e_1=e_4  & e_2 e_2=e_5 \\
{\bf A}_{131}&:&(\mathfrak{N}_2)_{5,11} &:& e_1 e_1 = e_5 & e_1 e_2=e_3 & e_2 e_1=e_4  & e_2 e_2=e_5 \\
{\bf A}_{132}&:&(\mathfrak{N}_2)_{5,12} &:& e_1 e_1 = e_3 & e_1 e_2=e_4 & e_2 e_1=e_4  & e_2 e_2=e_5 \\
{\bf A}_{133}(\lambda)&:&(\mathfrak{N}_2)_{5,13} (\lambda \in {\bf k}) &:& e_1 e_1 = e_3+\lambda e_5 & e_1 e_2=e_3 & e_2 e_1=e_4  & e_2 e_2=e_5 
\end{array}$$   
\newpage

\begin{center}
\label{tab4}
{\bf Table 4.} 
\end{center}
\medskip
{
$$\begin{array}{|l|l|l|l|}
\hline
\textrm{Bilinear pair} & \textrm{Bilinear map / algebra multiplication (g)} & \textrm{Equivalence group/ Automorphism group} \\
\hline
\hline
\begin{array}{l}{\bf A}_1 (\alpha), \alpha \in {\bf k}  \end{array}& 
\begin{array}{ll}
e_1 e_1 = e_1+e_2, & e_1 e_2=\alpha e_2,\\ 
e_2 e_2=0, & e_2 e_1=(1-\alpha)e_2 
\end{array} & 
\begin{array}{l}\left\{\begin{pmatrix} 1 & 0 \\ x & 1 \\ \end{pmatrix}: x\in {\bf k}\right\}  \end{array} \Tstrut\Bstrut\\
\hline
\begin{array}{l}{\bf A}_2   \end{array}& \begin{array}{ll} 
e_1 e_1 = e_2, & e_1 e_2= e_2,\\
e_2 e_2=0, & e_2 e_1=-e_2  \end{array} &
\begin{array}{l}
\left\{\begin{pmatrix} 1 & 0 \\ x & 1 \\ \end{pmatrix}: x\in {\bf k}\right\}  \end{array}\Tstrut\Bstrut \\
\hline
\begin{array}{l}{\bf A}_3   \end{array}
& \begin{array}{ll} 
e_1 e_1 = e_2, & e_1 e_2=0,\\
e_2 e_2=0, & e_2 e_1=0  \end{array} &
\begin{array}{l}\left\{\begin{pmatrix} x & 0 \\ y & x^2 \\ \end{pmatrix}: x\in {\bf k}^* \textrm{ and } y\in {\bf k}\right\}  \end{array}\Tstrut\Bstrut\\
\hline
\begin{array}{l}{\bf A}_4 (\alpha), \alpha \in {\bf k}_{\geq0}  \end{array}&
\begin{array}{ll} 
e_1 e_1 = \alpha e_1+e_2, & e_1 e_2= e_1+\alpha e_2,\\
e_2 e_2=0, & e_2 e_1=-e_1  \end{array}  &
\begin{array}{l}\Tstrut\left\{Id, \begin{pmatrix} -1 & 0 \\ 0 & 1 \\ \end{pmatrix} \right\} \mbox{ if }\alpha=0,\Bstrut\\
\,\left\{Id\right\} \mbox{ otherwise}\bstrut  \end{array}\\
\hline
\begin{array}{l}{\bf B}_1(\alpha), \alpha\in {\bf k}  \end{array}&
\begin{array}{ll} 
e_1 e_1 = 0, & e_1 e_2=(1- \alpha) e_1 + e_2,\\
e_2 e_2=0, & e_2 e_1=\alpha e_1- e_2  \end{array}  &
\begin{array}{l}\,\left\{Id\right\}  \end{array}\\
\hline
\begin{array}{l}{\bf B}_2(\alpha), \alpha \in {\bf k}  \end{array}&
\begin{array}{ll} 
e_1 e_1 = 0, & e_1 e_2=(1- \alpha) e_1,\\
 e_2 e_2=0, & e_2 e_1=\alpha e_1   \end{array} &
\begin{array}{l}\left\{\begin{pmatrix} x & 0 \\ 0 & 1 \\ \end{pmatrix}: x\in {\bf k}^*\right\}  \end{array}\Tstrut\Bstrut\\
\hline
\begin{array}{l}{\bf B}_3  \end{array}
& \begin{array}{ll} 
e_1 e_1 = 0, & e_1 e_2= e_2, \\
e_2 e_2=0, & e_2 e_1=-e_2  \end{array} &
\begin{array}{l}\left\{\begin{pmatrix} 1 & 0 \\ x & y \\ \end{pmatrix}: x\in {\bf k}, y\in {\bf k}^*\right\}  \end{array}\Tstrut\Bstrut \\
\hline 
\begin{array}{l}{\bf C}(\alpha, \beta),\\
(\alpha, \beta)\in {\bf k} \times {\bf k}_{\geq 0}  
\end{array}&
\begin{array}{ll} 
e_1 e_1 = e_2, & e_1 e_2=(1-\alpha)e_1+ \beta e_2, \\
e_2 e_2=e_2, & e_2 e_1=\alpha e_1 - \beta e_2 \end{array} &
\begin{array}{l}\Tstrut\left\{Id, \begin{pmatrix} -1 & 0 \\ 0 & 1 \\ \end{pmatrix} \right\}
\mbox{ if }\beta=0,\Bstrut\\
\,\left\{Id\right\}\mbox{ if }\beta\not=0\bstrut   \end{array}\\

\hline
\begin{array}{l}{\bf D}_1(\alpha, \beta),\\
(\alpha, \beta)\in \mathcal U  
\end{array}&
\begin{array}{ll} 
e_1 e_1 = e_1, & e_1 e_2=(1-\alpha)e_1+ \beta e_2,\\
e_2 e_2=0, & e_2 e_1=\alpha e_1 - \beta e_2  \end{array} &
\begin{array}{l}
\Tstrut\left\{Id, \begin{pmatrix} 1 & 1 \\ 0 & -1 \\ \end{pmatrix} \right\}
\mbox{ if }\beta=2\alpha-1, \Bstrut\\
\,\left\{Id\right\} \mbox{ if } \beta\not=2\alpha-1\bstrut 
\end{array}\\

\hline
\begin{array}{l}{\bf D}_2(\alpha, \beta),\\
(\alpha,\beta) \in {\bf k}^2\backslash \mathcal T   
\end{array}&
\begin{array}{ll} 
e_1 e_1 = e_1, & e_1 e_2=\alpha e_2,\\
e_2 e_2=0, & e_2 e_1= \beta e_2  \end{array} &
\begin{array}{l}\left\{\begin{pmatrix} 1 & 0 \\ 0  & x \\ \end{pmatrix}, x\in {\bf k}^*\right\}  \end{array}\Tstrut\Bstrut\\

\hline
\begin{array}{l}{\bf D}_3(\alpha, \beta),\\
(\alpha,\beta)\in {\bf k}^2\backslash \mathcal T 
\end{array}&
\begin{array}{ll}
e_1 e_1 = e_1, & e_1 e_2=e_1+ \alpha e_2,\\
e_2 e_2=0, & e_2 e_1=- e_1 + \beta e_2  \end{array} &
\begin{array}{l}\,\left\{Id\right\}  \end{array}\\

\hline
\begin{array}{l}{\bf E}_1(\alpha, \beta, \gamma, \delta),\\
(\alpha, \beta, \gamma, \delta) \in \mathcal V  
\end{array}&
\begin{array}{ll}
e_1 e_1 = e_1, & e_1 e_2=\alpha e_1+ \beta e_2,\\
e_2 e_2=e_2, & e_2 e_1=\gamma e_1 + \delta e_2  \end{array} &
\begin{array}{l}
\Tstrut\left\{Id, \begin{pmatrix} 0 & 1 \\ 1 & 0 \\ \end{pmatrix} \right\} \mbox{ if } (\alpha,\gamma)=(\delta,\beta)\neq(-1,-1),\Bstrut\\
\,S_3 \mbox{ if } (\alpha,\gamma)=(\delta,\beta)=(-1,-1),\bstrut\\
\,\left\{Id\right\} \mbox{ otherwise}\bstrut 
\end{array}\\

\hline
\begin{array}{l}{\bf E}_2(\alpha, \beta, \gamma), \\
(\alpha,\beta, \gamma)\in {\bf k}^3\backslash {\bf k} \times \mathcal T 
\end{array}
& \begin{array}{ll}e_1 e_1 = e_1, & e_1 e_2=(1-\alpha)e_1+ \beta e_2,\\
e_2 e_2=e_2,& e_2 e_1= \alpha e_1 + \gamma e_2  \end{array} &
\begin{array}{l}\,\left\{Id\right\}  \end{array}\\

\hline 
\begin{array}{l}{\bf E}_3(\alpha, \beta, \gamma),\\ 
(\alpha, \beta, \gamma)\in {\bf k}^2\times {\bf k}^*_{>1} 
\end{array}&
\begin{array}{ll}
e_1 e_1 = e_1, & e_1 e_2=(1-\alpha)\gamma e_1+ \frac{\beta}{\gamma} e_2,\\
e_2 e_2=e_2, & e_2 e_1=\alpha \gamma e_1 + \frac{1-\beta}{\gamma}e_2  \end{array}
& \begin{array}{l}\Tstrut\left\{Id, \begin{pmatrix} 0 & 1 \\ 1 & 0 \\ \end{pmatrix} \right\}
\mbox{ if } \gamma=-1 \mbox{ and }\alpha=\beta,\Bstrut\\
\,\left\{Id\right\} \mbox{ otherwise} \bstrut
\end{array}\\

\hline
\begin{array}{l}{\bf E}_4  \end{array}& 
\begin{array}{ll}
e_1 e_1 = e_1, & e_1 e_2=e_1+ e_2,\\
e_2 e_2=e_2, &  e_2 e_1=0   \end{array} 
& \begin{array}{l}
\,\left\{Id\right\} 
\end{array}\\

\hline
\begin{array}{l}{\bf E}_5(\alpha), \alpha \in {\bf k}  \end{array}&
\begin{array}{ll} 
e_1 e_1 = e_1, & e_1 e_2=(1-\alpha)e_1+ \alpha e_2,\\
e_2 e_2=e_2, & e_2 e_1=\alpha e_1+ (1-\alpha) e_2  \end{array} &
\begin{array}{l} \left\{\begin{pmatrix} x & y \\ 1-x & 1-y \\ \end{pmatrix}, x, y\in {\bf k}, x\not=y\right\} 
\end{array}\Tstrut\Bstrut\\
\hline 

\begin{array}{l}\mathfrak{N}_2  \end{array}&
\begin{array}{ll} 
e_1 e_1 = 0, & e_1 e_2=0,\\
e_2 e_2=0, & e_2 e_1=0  \end{array} & 
\,\,\,\,GL(2,{\bf k})\\ 

\hline
\end{array}$$
}

\newpage


\end{document}